\documentclass{article}

\usepackage{amssymb,latexsym,amsfonts}

\title{Darboux evaluations of algebraic Gauss
hypergeometric functions}

\author{Raimundas Vid\=unas\thanks{Supported by 
the Dutch NWO project 613-06-565, by the ESF NOG-project, 
and by the JSPS grant No 20740075.}\\
 \em Kobe University}

\newtheorem{theorem}{Theorem}[section]
\newtheorem{lemma}[theorem]{Lemma}

\newtheorem{definition}[theorem]{Definition}

\usepackage{latexsym}

\newcommand{\negvspace}{\vspace{-1pt}}

\newcommand{\hpg}[5]{{}_{#1}\mbox{\rm F}_{\!#2}\!
  \left(\left.{#3 \atop #4}\right| #5 \right) }

\newcommand{\hpgo}[2]{{}_{#1}\mbox{\rm F}_{\!#2}}

\newcommand{\proof}{{\bf Proof. }}
\newcommand{\qed}{\hfill $\Box$\\}
\newcommand{\equal}{\!\!\!=\!\!\!}
\newcommand{\CC}{{\Bbb C}}
\newcommand{\QQ}{{\Bbb Q}}
\newcommand{\PP}{{\Bbb P}}

\newcommand{\ZZ}{{\Bbb Z}}

\catcode`\@=11 \@addtoreset{equation}{section} \catcode`\@=12

\date{}

\begin{document}

\maketitle

\begin{abstract} 
Algebraic Gauss hypergeometric functions can be expressed explicitly
in several ways. One attractive way is to pull-back their hypergeometric equations
(with a finite monodromy) to Fuchsian equations with a finite cyclic monodromy, 
and express the algebraic  solutions as radical functions on the covering curve.
This article presents these pull-back transformations of minimal degree 
for the hypergeometric equations with the tetrahedral, octahedral
or icosahedral projective monodromy. The minimal degree is 4, 6 or 12, respectively.
The covering curves are called {\em Darboux curves}, and they have genus 0
or (for some icosahedral Schwarz types) genus 1.
\end{abstract}

\section{Introduction} \label{alghypers}

Algebraic Gauss hypergeometric functions were studied by many
authors; see for example \cite{schwarz72,fuchs,brioschi,
klein78,pepin,boulang,katz72,baltadw,singulm2,putulm}. This purpose of this
article is to derive and present satisfying explicit forms for these functions.

The Gauss hypergeometric equation is
\begin{equation} \label{hpgdee}
X\,(1-X)\,\frac{d^2Y(X)}{dX^2}+
\big(c-(a\!+\!b\!+\!1)\,Z\big)\,\frac{dY(X)}{dX}-a\,b\,Y(X)=0.
\end{equation}
This is a Fuchsian equation on $\PP^1_X$ with 3 regular singular points $X=0$,
$1$, $\infty$. The local exponent differences are:
$1-c$ at $X=0$; $c-a-b$ at $X=1$; and $a-b$ at $X=\infty$. 
Let $E(e_1,e_2,e_3)$ denote a hypergeometric equation
with the local exponent differences $e_1,e_2,e_3$ at the 3 singular points.
The order and multiplication by $-1$ of the exponent differences is disregarded
in this notation.

As well known \cite{schwarz72,baltadw,putulm}, 
hypergeometric equation (\ref{hpgdee}) has a basis of algebraic solutions if
and only if its monodromy group is finite. The following hypergeometric
equations (and their fractional-linear transformations; see Appendix
\ref{hpgdeqs}) have this property and are called {\em standard
hypergeometric equations} with algebraic solutions:
\begin{itemize}
\item $E(1,1/n,1/n)$, where $n$ is a positive integer. The hypergeometric
equation degenerates to a Fuchsian equation with two singular points. Its
monodromy group is a cyclic group with $n$ elements. \item
$E(1/2,1/2,1/n)$, where $n$ is an integer, $n\ge 2$. The projective 
monodromy group of this equation is the dihedral group with $2n$ elements.
\item $E(1/2,1/3,1/3)$. The projective monodromy 
is the tetrahedral group, isomorphic to $A_4$.
\item $E(1/2,1/3,1/4)$. The projective monodromy 
is the octahedral group, isomorphic to $S_4$.
\item $E(1/2,1/3,1/5)$. The projective monodromy 
is the icosahedral group, isomorphic to $A_5$.
\end{itemize}
Recall that the projective monodromy group is the monodromy group
(or the Galois group, if solutions are algebraic) of $y_1/y_2$,
where $y_1,y_2$ is a basis of solutions of the differential equation.
The quotient  $y_1/y_2$ is known as a {\em Schwarz map}.
In general, the two-dimensional monodromy representation of second order 
Fuchsian equations is a cyclic extension of the projective monodromy.

The celebrated theorem of Klein \cite{klein77} states that if a second order
linear homogeneous differential equation only has algebraic solutions, then
that equation is a {\em pull-back transformation} of a standard
hypergeometric equation from the list above. Explicitly, if the Fuchsian
equation has coefficients in $\CC(X)$, the pull-back transformation changes
the variable $Z$ in (\ref{hpgdee}) to a rational function $\varphi(X)$. 
In geometric terms, we have a finite covering $\varphi:\PP^1_X\to\PP^1_Z$
between two projective lines, and we {\em pull-back} the standard
hypergeometric equation from $\PP^1_Z$ onto $\PP^1_X$.

The possible finite projective monodromies for the second order Fuchsian equations
are the same as listed above: a cyclic, a dihedral, or the tetrahedral, octahedral 
or icosahedral groups. The standard hypergeometric equation in Klein's theorem
can be chosen to have the same projective monodromy. 
In particular, the theorem of Klein implies that if a hypergeometric equation
only has algebraic solutions, then it is a pull-back transformation of
a listed standard hypergeometric equation with the same projective monodromy.

Hypergeometric equations with finite monodromies were first classified
by Schwarz in \cite{schwarz72}. Disregarding hypergeometric equations with
a cyclic monodromy, Schwarz gave a list of 15 types of these
hypergeometric equations. One type consists of hypergeometric equations with
a dihedral projective monodromy. The other types are represented by the following
hypergeometric equations:
\begin{itemize}
\item $E(1/2,1/3,1/3)$, $E(1/3,1/3,2/3)$. The projective monodromy 
is   tetrahedral. 
\item $E(1/2,1/3,1/4)$, $E(2/3,1/4,1/4)$. The projective monodromy 
is octahedral. 
\item $E(1/2,1/3,1/5)$, $E(1/2,1/3,2/5)$,
$E(1/2,1/5,2/5)$, $E(1/3,1/3,2/5)$, $E(1/3,2/3,1/5)$, $E(2/3,1/5,1/5)$,
$E(1/3,2/5,3/5)$, $E(1/3,1/5,3/5)$, $E(1/5,1/5,4/5)$,
$E(2/5,2/5,2/5)$. The projective monodromy 
is icosahedral. 
\end{itemize}
We refer to Schwarz type of hypergeometric equations with algebraic
solutions by the triple of the exponent differences
$e_1,e_2,e_3$ 
of these representative equations. (Usually, the Schwarz type is denoted by a
roman numeral from I to XV.) We refer to the listed 14 hypergeometric
equations as {\em main representatives} of the Schwarz types.
Hypergeometric equations of the same Schwarz type are characterized by the
property that their hypergeometric solutions are {\em contiguous} to
hypergeometric solutions of the main representative (see Appendix \ref{contigrels}).

Algebraic solutions of differential equations can be represented in several ways.
For example, by minimal polynomial equations their satisfy \cite{singulm2}, 
or (if the Galois group is solvable) by nested radical expressions. 
In \cite{WeHuBe},  \cite{kleinvhw}  \cite[Chapter 1]{maintphd} an algorithm 
is developed to represent algebraic solutions of second order linear differential equations 
using Klein's theorem. The representation form is
\begin{equation} \label{kleinev}
\theta(x)\,H(\varphi(x)),
\end{equation}
where $\varphi(x)$ is a rational function, $\theta(x)$ is a radical function 
that define Klein's pull-back transformation, and $H(z)$ is a solution of a corresponding standard
hypergeometric equation.  The degree of $\varphi(x)$ is equal to the order of the projective monodromy
group, that is, 12, 24 or 60 for the three most interesting monodromies.

We propose to pull-back a hypergeometric equations with a finite monodromy
to Fuchsian equations with a cyclic monodromy group. Then some
hypergeometric solutions are transformed to rather simple radical functions.
We call a pull-back covering $\phi:D\to \PP^1$ of this kind a {\em Darboux
covering}. The covering curve $D$ is called a {\em Darboux curve}.
Identification of algebraic Gauss hypergeometric functions with radical
functions on an algebraic curve offers satisfying geometric intuition,
especially when the Darboux covering has low degree. We use the term {\em
Darboux evaluation} to refer to the aspired identification of hypergeometric
functions with radical functions.

The Darboux coverings for hypergeometric equations with finite dihedral projective
mondromies have degree 2. They are simply $Z=x^2$ if half-integer exponent 
differences are at $Z=0$ and $Z=\infty$. In \cite{tdihedral},
the pulled-back dihedral hypergeometric functions are expressed in terms
of power functions and terminating Appell's (double hypergeometric) $F_2$ or $F_3$ sums.

Theory of Darboux coverings is developed in \S \ref{darbouxcurves}. It
turns out that Darboux coverings for hypergeometric equations of the same
Schwarz type are identical. Therefore we have finitely many different
Darboux coverings and Darboux curves. We compute all Darboux coverings of
minimal degree, which turns out to be 4, 6 or 12 for the tetrahedral,
octahedral and icosahedral types, respectively. The corresponding Darboux
curves have genus 0 or (for some icosahedral types) genus 1. For each
Schwarz type, we use Darboux coverings of minimal degree and compute
Darboux evaluations for two hypergeometric solutions of the main
representative equation.
The evaluations are presented in \S \ref{database}.
Using these formulas, differentiation and contiguous relations, one can compute Darboux
evaluations for two different hypergeometric solutions of any hypergeometric
equation with the tetrahedral, octahedral or icosahedral projective monodromy.

Section \ref{databasec} presents the method used to compute Darboux  evaluations.
Most attention is paid to describing complicated computations on genus 1 Darboux curves. 
The computer algebra system {\sf Maple} was used in the computations.
Appendix \S \ref{appendix} gives short introductions  to pull-back transformations, 
Fuchsian equations, contiguous relations and other relevant topics.

\section{Hypergeometric formulas}
\label{database}

This section presents the main explicit results: Darboux evaluations for
two hypergeometric functions 
of each Schwarz type with the tetrahedral, octahedral and icosahedral projective monodromies. 
Each pair of hypergeometric  evaluations forms a solution basis of the main representative equation
of the respective Schwarz type, up to a power factor to the second evaluation.
Darboux evaluations for other hypergeometric functions of the same Schwarz type
can be computed by using contiguous relations if another evaluation
in the same contiguous orbit is known. Another contiguous evaluation
can be computed by differentiation.

We use Darboux coverings of minimal possible degree, which is 4, 6 and 12 
for the three considered monodromies, respectively, as shown in Lemma \ref{darbdegs}.
For the tetrahedral and octahedral Schwarz types, Darboux coverings are
evident from the arguments of the hypergeometric functions. Darboux
coverings for icosahedral Schwarz types are given in (\ref{isophi1}),
(\ref{isophi2}), (\ref{isophi3}), (\ref{isophi4}), (\ref{isophi5}),
(\ref{isophi6}). Some of these Darboux coverings are valid for two different
Schwarz types. For seven icosahedral types, the Darboux curve has genus 1
rather than 0. Weierstrass equations
$\xi^2=x\,(1+\alpha\,x+\beta\,x^2)$, 
with $\alpha,\beta\in\CC$, 
for these curves are given in (\ref{darbouxc1}), (\ref{darbouxc2}),
(\ref{darbouxc3}), (\ref{darbouxc4}).
The Darboux evaluations hold locally around $x=0$
or, if the Darboux curve has genus 1, around the point $(x,\xi)=(0,0)$. 
A simple way to check each evaluation is to expand both sides in power
series around $x=0$. If the Darboux curve has genus 1, one has to replace
$\xi$ by the respective $\sqrt{x}\,\sqrt{1+\alpha\,x+\beta\,x^2}$, and
expand the power series in $\sqrt{x}$.

\subsection{Tetrahedral hypergeometric equations}
\label{tetrahedral}

Here are two hypergeometric evaluations for solutions of  $E(1/2,1/3,1/3)$:
\begin{eqnarray} \label{fptetra1}
\hpg{2}{1}{1/4,-1/12}{2/3}{\frac{x\,(x+4)^3}{4(2x-1)^3}} 
& = & \big(1-2x\big)^{-1/4},\\ 
\label{fptetra1b}\hpg{2}{1}{1/4,\;7/12}{4/3}{\frac{x\,(x+4)^3}{4(2x-1)^3}}
& = & \frac{1}{1+\frac14x}\,\left(1-2x\right)^{3/4}.
\end{eqnarray}
If we multiply the second formula by $x^{1/3}(x+4)/(2x-1)$, 
the two evaluated functions form a basis of solutions of the same 
hypergeometric equation.
Here are two hypergeometric evaluations for solutions of  $E(1/3,1/3,2/3)$,
of the other tetrahedral Schwarz type:
\begin{eqnarray} \label{fptetra2}
\hpg{2}{1}{1/2,-1/6}{2/3}{\frac{x\,(x+2)^3}{(2x+1)^3}}
& = & \big(1+2x\big)^{-1/2},\\
\hpg{2}{1}{1/6,\;5/6}{4/3}{\frac{x\,(x+2)^3}{(2x+1)^3}} & = &
\frac{1}{1+\frac{1}{2}x}\,\big(1+2x\big)^{1/2}\,\big(1+x\big)^{1/3}. 
\end{eqnarray}

\subsection{Octahedral hypergeometric equations}

Here are two representative solutions of $E(1/2,1/3,1/4)$:
\begin{eqnarray} \label{fpocta1}
\hpg{2}{1}{7/24,-1/24}{3/4}{\frac{108\,x\,(x-1)^4}{(x^2+14x+1)^3}}
&=& \big(1+14x+x^2\big)^{-1/8},\\
\hpg{2}{1}{5/24,\;13/24}{5/4}{\frac{108\,x\,(x-1)^4}{(x^2+14x+1)^3}}
&=& \frac{1}{1-x}\,\big(1+14x+x^2\big)^{5/8}.
\end{eqnarray}
And here are two representative solutions of $E(1/4,1/4,2/3)$:
\begin{eqnarray} \label{fpocta2}
\hpg{2}{1}{7/12,-1/12}{3/4}{\frac{27\,x\,(x+1)^4}{2(x^2+4x+1)^3}}
& = & \frac{\left(1+\frac12x\right)^{1/4}}{(1+4x+x^2)^{1/4}},\\
\hpg{2}{1}{1/6,\;5/6}{5/4}{\frac{27\,x\,(x+1)^4}{2(x^2+4x+1)^3}}
& = & \frac{\big(1+2x\big)^{1/4}\,\big(1+4x+x^2\big)^{1/2}}{1+x}.
\end{eqnarray}

\subsection{Icosahedral hypergeometric equations}

The Darboux covering for hypergeometric equations of the Schwarz types
$(1/2,1/3,1/5)$ and $(1/2,1/3,2/5)$ is
\begin{equation} \label{isophi1}
\varphi_1(x)=\frac{1728\;x\;(x^2-11x-1)^5}{(x^4+228x^3+494x^2-228x+1)^3}
\end{equation}
Representative evaluations for solutions of $E(1/2,1/3,1/5)$ are:
\begin{eqnarray}
\hpg{2}{1}{19/60,-1/60}{4/5}{\varphi_1(x)} & \equal &
\big(1-228x+494x^2+228x^3+x^4\big)^{-1/20},\\
\hpg{2}{1}{11/60,\,31/60}{6/5}{\varphi_1(x)} & \equal &
\frac{(1\!-\!228x\!+\!494x^2\!+\!228x^3\!+\!x^4)^{11/20}}{1+11x-x^2}.
\end{eqnarray}
Similarly, basic evaluations for solutions of $E(1/2,1/3,2/5)$ are:
\begin{eqnarray}
\hpg{2}{1}{13/60,-7/60}{3/5}{\varphi_1(x)} & \equal &
\frac{1-7x}{\big(1-228x+494x^2+228x^3+x^4\big)^{7/20}},\\
\hpg{2}{1}{17/60,\,37/60}{7/5}{\varphi_1(x)} & \equal &
\frac{\left(1+\frac{1}{7}x\right)(1\!-\!228x\!+\!494x^2\!+\!228x^3\!+\!x^4)^{17/20}}
{\left(1+11x-x^2\right)^2}.\quad
\end{eqnarray}

The Darboux covering for Schwarz type $(1/2,1/5,2/5)$ is the following:
\begin{equation} \label{isophi2}
\varphi_2(x)= \frac{64\,x\,(x^2-x-1)^5}{(x^2-1)\,(x^2+4x-1)^5}.
\end{equation}
Representative evaluations for solutions of $E(1/2,1/5,2/5)$ are:
\begin{eqnarray}
\hpg{2}{1}{7/20,-1/20}{4/5}{\varphi_2(x)} & \equal &
\frac{(1+x)^{7/20}}{(1-x)^{1/20}\,(1-4x-x^2)^{1/4}},\\
\hpg{2}{1}{3/20,\,11/20}{6/5}{\varphi_2(x)} & \equal &
\frac{(1+x)^{3/20}\,(1-x)^{11/20}\,(1-4x-x^2)^{3/4}}{1+x-x^2}. 
\end{eqnarray}

Darboux curves for other icosahedral Schwarz types have genus 1.
The Darboux curve for hypergeometric equations of the Schwarz types
$(1/3,1/3,2/5)$ and $(1/3,2/3,1/5)$ is given by the equation
\begin{equation} \label{darbouxc1}
C_3: \qquad \xi^2=x\;(1+33x-9x^2).
\end{equation}
The Darboux covering is
\begin{equation} \label{isophi3}
\varphi_3(x,\xi)=\frac{144\;\xi\,\left(1+33x-9x^2\right)^2\,(1-9\xi+54x)}
{\left(1+21\xi-117x+9x\xi-234x^2\right)^3}
\end{equation}
Here are basic evaluations for solutions of $E(1/3,1/3,2/5)$:
\begin{eqnarray} \label{icosellipta}
\hpg{2}{1}{3/10,-1/30}{3/5}{\varphi_3(x,\xi)} & \equal &
\frac{(1-9\xi+54x)^{1/30}}{(1+21\xi-117x+9x\xi-234x^2)^{1/10}},\\
\hpg{2}{1}{7/10,\;11/30}{7/5}{\varphi_3(x,\xi)} & \equal &
\!\frac{\left(1+21\xi-117x+9x\xi-234x^2\right)^{11/10}}
{(1-9\xi+54x)^{11/30}\,(1+33x-9x^2)}. \label{icoselliptc3}
\end{eqnarray}
And here are basic evaluations for solutions of $E(1/3,2/3,1/5)$:
\begin{eqnarray} \label{icoselliptk3}
\hpg{2}{1}{\!17/30,-1/10}{4/5}{\varphi_3(x,\xi)}\! & \equal & \frac{
(1-9\xi+54x)^{13/30}}{\left(1+21\xi-117x+9x\xi-234x^2\right)^{3/10}},\\  \label{icoselliptz3}
\hpg{2}{1}{1/10,\,23/30}{6/5}{\varphi_3(x,\xi)}\! & \equal & \!\frac{
\left(1\!+\!21\xi\!-\!117x\!+\!9x\xi\!-\!234x^2\right)^{3/10}\!
(1\!-\!9\xi\!+\!54x)^{7/30}(\xi\!+\!5x)}{\xi\;(1+9x)}.\label{icoselliptm3}
\nonumber\\ 
\end{eqnarray}

The Darboux curve for the Schwarz types $(2/3,1/5,1/5)$ and $(1/3,2/5,3/5)$
is given by
\begin{equation} \label{darbouxc2}
C_4: \qquad \xi^2=x\;(1+5x-5x^2).
\end{equation}
The Darboux covering is 
\begin{equation} \label{isophi4}
\varphi_4(x,\xi)=\frac{432\,x\,\left(1-\frac{7}{5}\xi-9x-x^2\right)^5\,
(1\!+\!50x\!-\!125\xi^2\!+\!450x\xi\!-\!500x^2)}{(5\xi\!+\!57x)\,
\left(1\!+\!\frac{18}{5}\xi\!-\!16x\!+\!x^2\right)^5
(1\!+\!50x\!-\!125\xi^2\!-\!450x\xi\!-\!500x^2)}.
\end{equation}
Basic evaluations for solutions of $E(2/3,1/5,1/5)$ are:
\begin{eqnarray} \label{icosellipta4}
\hpg{2}{1}{1/6,-1/30}{4/5}{\varphi_4(x,\xi)}\! & \equal & \!\frac{
\left(1-\frac35\xi-\frac{34}{5}x\right)^{1/6}}
{(1\!+\!3\xi\!-\!20x)^{1/6}\left(
1\!+\!50x\!-\!125\xi^2\!-\!450x\xi\!-\!500x^2\right)^{1/30}},
\hspace{25pt} \\  \label{icoselliptc4}
\hpg{2}{1}{1/6,\,11/30}{6/5}{\varphi_4(x,\xi)}\! & \equal &
\!\left(1\!+\!50x\!-\!125\xi^2\!-\!450x\xi\!-\!500x^2\right)^{11/30} \times
\nonumber \\ && \! 
\frac{(1+3\xi-20x)^{5/6}\left(1-\frac35\xi-\frac{34}5x\right)^{1/6}
\left(1+\frac{21}4\xi+\frac{41}4x\right)}
{(1-9x)\,\left(1-\frac74\xi-\frac{15}2x\right)\,(1+5\xi+10x)}.
\end{eqnarray}
Basic evaluations for solutions of $E(1/3,2/5,3/5)$ are:
\begin{eqnarray}\label{icoselliptk4}
\hpg{2}{1}{\! 13/30,-1/6}{3/5}{\varphi_4(x,\xi)}\! & \equal & \!\frac{
(1\!+\!50x\!-\!125\xi^2\!-\!450x\xi\!-\!500x^2)^{13/30}
\left(1-3\xi+2x\right)}{\left(1+3\xi-20x\right)^{5/6}\,
\left(1-\frac35\xi-\frac{34}{5}x\right)^{1/6}\,(1+5\xi+10x)}, 
\hspace{25pt} \\  \label{icoselliptz4}
\hpg{2}{1}{5/6,\,7/30}{7/5}{\varphi_4(x,\xi)}\! & \equal &
\left(1\!+\!50x\!-\!125\xi^2\!-\!450x\xi\!-\!500x^2\right)^{7/30} \times
\nonumber \\ && \hspace{0pt} \label{icoselliptm4}
\frac{\left(1\!+\!\frac{18}{5}\xi\!-\!16x\!+\!x^2\right)^{7/6}
\!\left(1\!+\!\frac{1}{25}x\right)^{5/6}\!(1\!+\!5\xi\!+\!10x)}
{\left(1-\frac{7}{5}\xi-9x-x^2\right)^2\;(1-5x)^{7/6}}.
\end{eqnarray}

The Darboux curve for the Schwarz type $(1/3,1/5,3/5)$ is given by
\begin{equation} \label{darbouxc3}
C_5: \qquad \xi^2=x\;(1+x)\,(1+16x).
\end{equation}
The Darboux covering is
\begin{equation} \label{isophi5}
\varphi_5(x,\xi)=-\frac{54\;(\xi+5x)^3\,(1-2\xi+6x)^5}
{(1\!-\!16x^2)\,(\xi-5x)^2\,(1-2\xi-14x)^5}.
\end{equation}
Representative evaluations for solutions of $E(1/3,1/5,3/5)$ are:
\begin{eqnarray} \label{icosellipta5}
\hpg{2}{1}{8/15,-1/15}{4/5}{\varphi_5(x,\xi)}\! & \equal &
\frac{(1+4x)^{8/15}\,(\xi+5x)^{1/6}\;x^{1/15}}
{(1-2\xi-14x)^{1/3}\;(\xi-3x)^{3/10}},\\  \label{icoselliptz5}
\hpg{2}{1}{2/15,\,11/15}{6/5}{\varphi_5(x,\xi)}\! & \equal & \frac{
(1\!-\!\xi\!+\!x)(1\!-\!2\xi\!-\!14x)^{2/3}(\xi\!+\!5x)^{1/6}(\xi\!-\!3x)^{13/10}}
{(1+\xi+x)\,(1-2\xi+6x)\,(1+4x)^{13/15}\;x^{11/15}}. \hspace{20pt}
\end{eqnarray}

The Darboux curve for the Schwarz types $(1/5,1/5,4/5)$ and  $(2/5,2/5,2/5)$
is given by
\begin{equation} \label{darbouxc4}
C_6: \qquad \xi^2=x\;(1+x-x^2).
\end{equation}
The Darboux covering is 
\begin{equation} \label{isophi6}
\varphi_6(x,\xi)=\frac{16\;\xi\,\left(1+x-x^2\right)^2\,(1-\xi)^2}
{(1+\xi+2x)\,(1+\xi-2x)^5}.
\end{equation}
Basic evaluations for solutions of $E(1/5,1/5,4/5)$ are:
\begin{eqnarray} \label{icosellipta6}
\hpg{2}{1}{7/10,-1/10}{4/5}{\varphi_6(x,\xi)} & \equal & \frac{
(1-\xi+2x)^{1/15}\,(1-\xi)^{3/5}}{(1+\xi+2x)^{7/30}\,\sqrt{1+\xi-2x}},\\
\hpg{2}{1}{1/10,\,9/10}{6/5}{\varphi_6(x,\xi)} & \equal & \frac{
(\xi+2x+x^2)\;(1+\xi)^{1/10}\,(1-\xi)^{3/10}}{\xi\,(1-\xi+2x)^{1/30}
(1+\xi+2x)^{2/15}\sqrt{1+\xi-2x}}.\quad
\end{eqnarray}
Basic evaluations for solutions of $E(2/5,2/5,2/5)$ are:
\begin{eqnarray}
\hpg{2}{1}{3/10,-1/10}{3/5}{\varphi_6(x,\xi)}\! & \equal & \frac{
(1-\xi+2x)^{2/15}\,(1+\xi+2x)^{1/30}\,(1-\xi)^{1/5}}{\sqrt{1+\xi-2x}},\\  \label{icoselliptz}
\hpg{2}{1}{3/10,\,7/10}{7/5}{\varphi_6(x,\xi)} & \equal & \frac{
(1+\xi+2x)^{7/30}\,(1+\xi)^{1/5}\,(1+\xi-2x)^{3/2}}
{(1-\xi+2x)^{1/15}\,(1-\xi)^{2/5}\,(1+x-x^2)}.
\end{eqnarray}

\section{Darboux curves} \label{darbouxcurves}

In the context of integration theory of vector fields \cite{darbouxc}, \cite{moulino},
one considers a polynomial vector field on
$\CC^2$ given by a derivation ${\cal L}=f(x,y)\,\partial/\partial x
+g(x,y)\,\partial/\partial y$ with $f,g\in\CC[x,y]$. A polynomial $p(x,y)$
is called a {\em Darboux polynomial} for the vector field (or the derivation
$\cal L$) if $p(x,y)$ divides ${\cal L}p(x,y)$ in the ring $\CC[x,y]$. An
algebraic curve defined as the zero set of a Darboux polynomial is called a
{\em Darboux curve}. A Darboux curve is infinitesimally invariant under the
vector field. Hence an alternative term is {\em invariant algebraic curve},
as in \cite{darbouxphys} for example.

In differential Galois theory we have the following definition
\cite{jaweilphd}, \cite{Singer92}. Let $K$ be a differential field, and let
$R=K[y_1,\ldots,y_n]$ be a differential ring. Let $\cal D$ denote the
derivation on $R$, and suppose that it extends the derivation on $K$. Then
$p\in R$ is a {\em Darboux polynomial} for $\cal D$ if $p$ divides ${\cal
D}p$ in $R$. For example, consider differential equation
\begin{equation} \label{genlhde}
y^{(n)}+a_{n-1}\,y^{(n-1)}+\ldots+a_1\,y'+a_0\,y=0, \qquad \mbox{with} \quad
a_0,a_1,\ldots,a_{n-1}\in K.
\end{equation}
Let $\cal D$ be the derivation on $K[y,y',\ldots,y^{(n-1)}]$ defined by 
${\cal D} y=y'$, ${\cal D}y'=y'',\ldots$, ${\cal D} y^{(n-2)}=y^{(n-1)}$,
${\cal D} y^{(n-1)}=-a_{n-1}y^{(n-1)}-\ldots-a_1y'-a_0y$.
Darboux polynomials for this derivation correspond to 
semi-invariants of the differential Galois group for (\ref{genlhde});
see \cite[Theorem 38]{jaweilphd}. 
If the order $n$ of (\ref{genlhde}) is equal to $2$, then one considers
Darboux polynomials in $K[u]$ for the derivation $\cal D$ defined by ${\cal
D}u=-u^2-a_1u-a_0$. In \cite{jaweil94,UlmeWeil96} these Darboux
polynomials are called {\em special polynomials}.

As we see, the terms ``Darboux polynomials", ``Darboux curves"
are not consistently used though there is a common theme. 
In the context of algebraic solutions of differential equations,
we wish to use the term ``Darboux curves" to offer a geometric guidance. 
We offer the following definition of Darboux curves and Darboux coverings.
\begin{definition} \rm
Let $C$ denote an algebraic curve (see Appendix \ref{secalgcurves}). We
suppose that the function field $\CC(C)$ is a differential field. Consider
differential equation (\ref{genlhde}) on $\PP^1$, assuming $K=\CC(\PP^1)$.
We say that a finite covering $\phi:C\to\PP^1$ is a {\em Darboux covering}
for (\ref{genlhde}), if a pull-back of (\ref{genlhde}) with respect to
$\phi$ has a solution $Y$ such that:
\begin{enumerate}
\item Its logarithmic derivative $u=Y'/Y$ is in $\CC(C)$;%
\item The algebraic degree of $u$ over $K$ is precisely the degree of $\phi$.%
\end{enumerate}
In the described situation, $C$ is called a {\em Darboux curve} for
(\ref{genlhde}).
\end{definition}

To see connection with previous definitions, let us assume that the order
$n$ of (\ref{genlhde}) is equal to 2. We take the differential field $K=\CC(z)$ with $z'=1$.
Suppose that we have a Darboux covering $\phi:C\to\PP^1$ of degree $d$,
and let $u$ be a required logarithmic derivative in $\CC(C)$.
Then we have the following facts:
\begin{itemize}
\item The logarithmic derivative $u$ is an algebraic solution of the
associated Riccati equation $u'+u^2+a_1u+a_0=0$. Let $P(u)=0$ denote the
minimal monic polynomial equation defining $u$ over $K$. The polynomial $P$
has degree $d$, and it is a defining equation for the Darboux curve $C$.
Therefore $\CC(C)=\CC(z,u)$.
\item The polynomial $P$ is a Darboux polynomial for the 
derivation ${\cal D}u=-u^2-a_1u-a_0$ in the definition of differential Galois theory; 
see \cite[Lemma 2.4]{UlmeWeil96}.%
\item The expression $y^{n+1}P(y'/y)$ is a homogeneous polynomial in
$K[y,y']$. It is a Darboux polynomial for the corresponding derivation on
$K[y,y']$.%
\item Assume that $a_1,a_0\in\CC[z]$, and consider the vector field
$\partial/\partial z-(u^2+a_1u+a_0)\,\partial/\partial u$. Then $P$ is a
Darboux polynomial according the first definition above.
\end{itemize}
Besides, we have the following facts.
\begin{itemize}
\item Suppose that $P=u^m+\sum_{j=0}^{d-1} b_j\,u^j$.
Then $-b_{d-1}$ is the logarithmic derivative of an exponential solution of
the $d$-th symmetric tensor power of differential equation (\ref{genlhde}); 
see \cite[Theorem 2.1]{UlmeWeil96} or \cite[\S 3.2]{singulm2}. 
All other coefficients $b_j$ are determined by $b_{d-1}$ and 
the differential equation \cite[\S 2]{UlmeWeil96}.
\item The mentioned exponential solution (which can be expressed as
$\exp\int -b_{d-1}$) is a degree $d$ semi-invariant of the
differential Galois group \cite[\S 2]{UlmeWeil96}.%
\item Inside a Picard-Vessiot extension for (\ref{genlhde}), the field
$\CC(C)$ is fixed by a 1-reducible subgroup of the differential Galois group
for (\ref{genlhde}) of finite index \cite[Lemma 3.1]{singulm2}.%
\item If the differential Galois group of (\ref{genlhde}) is finite, the
1-reducible subgroups are cyclic subgroups \cite[Lemma 1.5]{UlmeWeil96}.
\end{itemize}
As we see, Darboux coverings correspond to algebraic Riccati solutions and 
semi-invariants of the differential Galois group. 

We defined Darboux curves and coverings for general linear differential equations. 
They can be computed from minimal polynomials for Riccati solutions or by direct
solution of pull-backed equations. Pulled-back equations on Darboux curves 
have cyclic monodromy groups, hence their solutions can be expressed in
terms of radical functions.  Solutions of the original equation can be presented
as those radical solutions of the transformed equation, keeping in mind the inverse 
substitution defined by the Darboux covering. We refer to these radical expressions 
as  {\em Darboux evaluations}  of solutions of the original equation.
They could be satisfying for many purposes.
For convenience, Darboux curves of genus 0 can be parametrized,
and Weierstrass forms of genus 1 Darboux used. 
Thanks to contiguous relations, Darboux evaluations for algebraic hypergeometric 
functions can be derived using only a few representative Darboux evaluations
for each Schwarz type.

Darboux evaluations for solutions of hypergeometric equations 
with dihedral projective monodromies 
is recognized in  \cite{tdihedral} or \cite[2.5.5]{bateman}. 
The Darboux coverings in the dihedral case are simple quadratic coverings,
such as $\varphi(x)=x^2$.
The Darboux curves and coverings for hypergeometric equations
were introduced earlier  in \cite[\S 4.2]{myphd}. 
The Darboux coverings of minimal degree for all tetrahedral, octahedral types
and for the icosahedral types $(1/2,1/3,1/5)$, $(1/2,1/3,2/5)$, $(1/2,1/5,2/5)$, $(1/3,1/3,2/5)$,
$(1/3,1/3,2/5)$ were computed in that PhD thesis.

\subsection{Basic properties of Darboux curves}

The following lemmas identify Darboux coverings for hypergeometric equations
of the same Schwarz type, and imply that there are only finitely many different
Darboux coverings for all hypergeometric equations. These are key attractive 
features for using Darboux pull-backs.
\begin{lemma} \label{dschwiso}
Suppose that $\phi:D\to\PP^1$ is a Darboux covering for a hypergeometric equation $E_1$
with a tetrahedral, octahedral or icosahedral monodromy, and $E_2$ is hypergeometric equation
of the same Schwarz type. Then $\phi:D\to\PP^1$ is a Darboux covering for $E_2$ as well.
\end{lemma}
\proof Let $y_1$ be a hypergeometric solution of $E_1$. Since $E_1$ and
$E_2$ have the same Schwarz type, there is a hypergeometric solution $y_2$
which is contiguous to $y_1$. Recall that $y_1'$ is contiguous to $y_1$ as well.
Therefore there is a contiguous relation $y_2=ay_1+by_1'$, with
$a,b\in\CC(z)$. There is also a contiguity relation $y_2'=cy_1+dy_1'$, with
$c,d\in\CC(z)$. Let $u_1$ denote the Riccati solution $y_1'/y_1$ for $E_1$,
and let $u_2$ denote the Riccati solution $y_2'/y_2$ for $E_2$. Then
$u_2=(c+du_1)/(a+bu_1)$. This implies that the function fields $\CC(z,u_1)$
and $\CC(z,u_2)$ are isomorphic. \qed

The following lemma categorizes all possible Darboux coverings for hypergeometric
equations with the considered monodromies.
\begin{lemma} \label{darbdegs}
Suppose that the projective monodromy 
$G$ of a hypergeometric equation $(\ref{hpgdee})$ is
either tetrahedral $A_4$, or octahedral $S_4$, or icosahedral $A_5$.
\begin{itemize}
\item If $G\cong A_4$ (tetrahedral group), the Darboux coverings 
have degree $4$, $6$ or $12$. 
\item If $G\cong S_4$ (octahedral group), the Darboux coverings 
have degree $6$, $8$, $12$ or $24$.
\item If $G\cong A_5$ (icosahedral group), the Darboux coverings 
have degree $12$, $20$, $30$ or $60$. 
\item The Darboux coverings of each degree are unique up to automorphisms 
of the Darboux curve. 
\end{itemize}
\end{lemma}
\proof  As noted, Darboux curves correspond to algebraic Riccati solutions.
According to \cite[Corollary 1.7]{UlmeWeil96}, the given numbers are possible
degrees of minimal polynomials for solutions of the Riccati equation
associated to the hypergeometric equations with considered monodromy groups.
For each projective monodromy group, Riccati solutions of the maximal degree 
are not unique, but they define the same function field (because of irreducibility
of the monodromy) and hence the same Darboux curve.
Irreducible polynomials for Riccati solutions of other (non-maximal) degrees 
are unique except that there are two irreducible polynomials of degree $4$ in the tetrahedral case.
To show that those two degree 4 polynomials define isomorphic Darboux curves,
we first use Lemma \ref{dschwiso} to reduce consideration to the two equations $E(1/2,1/3,1/3)$
and $E(1/3,1/3,2/3)$, representing the two distinct tetrahedral Schwarz types. 
The two points with the exponent differences $1/3$ of either of these equations are interchanged
by a fractional-linear transformation (say, $X\mapsto 1-X$), and the two minimal polynomials
differ by the same transformation.
\qed

The uniqueness claim does not mean that coverings for pull-backs of considered 
hypergeometric equations to Fuchsian equations with the isomorphic cyclic monodromies 
have to be isomorphic. In fact, Litcanu \cite[\S 2]{Litcanu2004} observed that there are two
pull-back transformations of $E(1/2,1/3,1/4)$ 
to a Fuchsian equation with the $\ZZ/2\ZZ$ monodromy.  
One pull-back covering is the degree 12 Darboux covering 
\begin{equation} \label{eq:octa2c2}
\frac{27(x-1)^4(x^2+6x+1)^4}{(x^2-10x+1)^3\,(3x^2+2x+3)^3},
\end{equation}
and the other is $108x^2(x^2-1)^4/(x^4+14x^2+1)^3$. 
The latter covering is a composition of the degree 6 covering for octahedral equations
(for a pull-back to $\ZZ/4\ZZ$) and a quadratic covering. It is not a Darboux covering,
because the respective Riccati solution is still of degree 6. The two coverings of degree 12
have different branching patterns.

\subsection{Darboux coverings for standard hypergeometric equations}
\label{darboux4st}

Here are the main facts about Darboux coverings for standard hypergeometric
equations.
\begin{lemma} \label{dramifico}
Let $H$ denote a standard hypergeometric equation with tetrahedral,
octahedral or icosahedral differential Galois group $G$. Let
$\varphi:D\to\PP^1$ be a Darboux covering for $H$ of degree $m$. Then:
\begin{enumerate}
\item[(i)] The Darboux curve $D$ has genus zero.
\item[(ii)] Let $\widetilde{\varphi}:\widetilde{D}\to\PP^1$ be a Darboux
covering for $H$ of the maximal degree $|G|$. Then we have a factorization
$\widetilde{\varphi}=\gamma\circ\varphi$, where the covering
$\gamma:\widetilde{D}\to D $ is given (up to fractional-linear
transformations of $D$ and $\widetilde{D}$) by $x\mapsto x^{|G|/m}$.
\item[(iii)] Let $X\in\PP^1$ be a regular singular point of the
hypergeometric equation. Assume that the local exponent difference at $X$
has denominator $k$. Then there are $\lfloor m/k\rfloor$ points above $X$
with branching index $k$, and other points above $X$ are unramified.
\end{enumerate}
\end{lemma}
\proof First we prove this theorem for the Darboux covering
$\widetilde{\varphi}$ of maximal degree $m=|G|$. The first statement follows
classically from \cite{icosaklein}. In particular, there is an action of $G$
on $\widetilde{D}\cong\PP^1$, and the projection $\PP^1\to\PP^1/G\cong\PP^1$
is precisely $\widetilde{\varphi}$. This map is also the inverse of a
Schwarz map for the hypergeometric equation. At each fiber of
$\widetilde{\varphi}$ the points have the same branching index, since
$\CC(\widetilde{D})\supset\CC[z]$ is a Galois extension. This extension is
also the Piccard-Vessiot extension for $H$, so a suitable pull-back of $H$
with respect to $\widetilde{\varphi}$ has trivial Galois group. Hence all
exponent differences of any pull-back of $H$ with respect to
$\widetilde{\varphi}$ are integers. As a consequence, in the situation of
part {\em (iii)} the branching indices should be integer multiples of $k$.
Since $g(\widetilde{D})=0$, Hurwitz theorem leaves only one possibility
which is described in part {\em (iii)}. 
Part {\em (ii)} is trivial in the considered case.

Now we consider a general Darboux covering $\varphi:D\to\PP^1$. Let $u$ be a
corresponding Riccati solution of degree $m$, so that $\CC(D)\cong\CC(z,u)$.
Let $K\supset\CC(z)$ denote the the Piccard-Vessiot extension of
(\ref{hpgdee}). We have $K\cong\CC(\widetilde{D})$ as just above. Consider
the tower of field extensions $K\supset\CC(D)\supset\CC(z)$. Here
$K\supset\CC(D)$ is the Piccard-Vessiot extension for the differential
equation $y'=uy$, so its Galois group must be a cyclic subgroup of $G$ of
index $m$. The existence of the corresponding covering
$\gamma:\widetilde{D}\to D$ implies {\em (i)}. Hurwitz theorem implies that
there are exactly two branching points of $\gamma$. This implies {\em (ii)}.
We see that the equation $y'=uy$ has exactly two singular points, so in the
situation of part {\em (iii)} branching indices are equal to either $k$ or
$1$, and that there in total only two unramified points above the three
regular singular points of (\ref{hpgdee}). This information gives part {\em
(iii)}. \qed

The Darboux coverings for the standard hypergeometric equations
can be computed from scratch, by using the branching pattern prescribed by
part {\em (iii)} of Lemma \ref{dramifico}, and Algorithm 1 in
\cite[\S 3]{thyperbolic} for example. On the other hand, the coefficients
in their rational expressions are familiar from the (semi)-invariants 
of the classical action of $A_4$, $S_4$, $A_5$ on $\CC[x]$ (or on
homogeneous polynomials in $\CC[x,y]$); see \cite{icosaklein}.

Explicit expressions for the Darboux coverings of the minimal degree 4, 6 or 12
are evident in formulas (\ref{fptetra1}), (\ref{fpocta1}), (\ref{isophi1}), respectively.
The covering of degree 6 for the standard tetrahedral equation is given by
\begin{equation}  \label{eq:tetra6}
X\mapsto \frac{\left(x^2-6\,x-3\right)^3}
{\left(x^2+6\,x-3\right)^3}.
\end{equation}
For the standard octahedral equation, the Darboux covering of degree 8 is given by
\begin{equation} \label{eq:octa8}
X\mapsto \frac{(x^2+20x-8)^4}{256\,x\,(x+1)^3\,(x-8)^3},
\end{equation}
while the covering of degree 12 is given in (\ref{eq:octa2c2}).
The Darboux coverings of degree 20 and 30 for the standard 
icosahedral equation are given by the rational functions
\begin{eqnarray} \label{icosacovs}
\frac{64\,(x^4+55x^3-165x^2-275x+25)^5}
{125x(x^2\!+\!5x\!+\!40)^3(x^2\!-\!40x\!-\!5)^3(8x^2\!-\!5x\!+\!5)^3},\nonumber\\
\frac{27\,(x^2+2x+5)^5\,(x^4+20x^3-210x^2+100x+25)^5}
{(3x^2\!-\!10x\!+\!15)^3(x^4\!+\!70x^2\!+\!25)^3
(x^4\!-\!60x^3\!-\!370x^2\!-\!300x\!+\!25)^3} .
\end{eqnarray}
The Darboux coverings of the maximal degree $D\in\{12,24,60\}$ 
can be obtained by composing an indicated Darboux covering of smaller degree $m$
with the cyclic covering $x\mapsto x^{D/m}$. Different rational expressions for
the maximal Darboux coverings are related to each other by M\"obius transformations,
though the coefficients are cumbersome algebraic numbers.
For example, to transform $x^3(x^3+4)^3/4(2x^3-1)^3$ 
to the composition of (\ref{eq:tetra6}) with $x\mapsto x^2$, 
a M\"obius transformation is 
\[
x\mapsto \sqrt[3]{5-3\sqrt3} \; \frac{x+\sqrt{3+2\sqrt3}}{x-\sqrt{-3+2\sqrt3}}.
\]

\subsection{Darboux curves for other Schwarz types}

The following lemma allows us to compute all necessary Darboux coverings
once Darboux coverings for standard hypergeometric equations (from \S \ref{darboux4st}) 
are known. For this purpose we usually take
$E_1$ to be a main representative equation of a chosen Schwarz type, and we
usually take $E_0$ to be the corresponding standard equation.
\begin{lemma} \label{dfiberprod}
Let $E_0$, $E_1$ denote two linear differential equations. Suppose that
$E_1$ is a pull-back of $E_0$ with respect to a covering
$\psi:\PP^1\to\PP^1$. Suppose that $\phi_0:D_0\to\PP^1$ is a Darboux
covering for $E_0$. Then the fiber product $D_1$ of $\psi:\PP^1\to\PP^1$ and
$\phi_0:D_0\to\PP^1$ is a Darboux curve for $E_1$, and the projection
$\phi_1:D_1\to\PP^1$ (to the source of $\psi$) is a Darboux covering for
$E_1$ of the same degree as $\phi_0$.
\end{lemma}
\proof Let $z$ denote a projective parameter for the $\PP^1$ below, and let
$x$ denote a projective parameter of the other $\PP^1$, so that the covering
$\psi$ corresponds to the extension $\CC(x)\supset\CC(z)$. Suppose that
$u_0$ be a Riccati solution for $E_0$ which determines $D_0$, so that
$\CC(D_0)\cong\CC(z,u_0)$. Then there is a solution $y$ for $E_0$ which
satisfies $y'=u_0y$. The logarithmic derivative of the pull-back of $y$ with
respect to (\ref{algtransf}) is $u_0\psi'+\theta'/\theta$. It is a Riccati
solution for $E_1$, and it lies in
$\CC(x,u_0)\cong\CC(x)\otimes_{\scriptsize\CC(z)}\CC(x,u_0)$ which is the
function field for the fiber product of $\psi$ and $\phi_0$. Since the
degree of the projection to $\PP^1$ is equal to the degree of $\psi$, the
claim follows.\qed

The simplest coverings for Klein pull-backs $\psi$ from standard
hypergeometric equations to other main representatives were first computed
in \cite{schwarz72}. These Klein coverings are familiar from classical
transformations of hypergeometric series (see for instance \cite{specfaar}
or \cite[\S 5]{algtgauss}), except the ones for the Schwarz types
$(1/2,1/3,2/5)$, $(1/3,2/3,1/5)$, $(1/3,2/5,3/5)$. Here are rational
functions that define Klein pull-backs for their main representatives:
\begin{equation} \label{okleinms}
\frac{x^2\,(189-64x)^5}{(3584x^2\!+\!2457x\!-\!2916)^3},\
\frac{4\,x\,(25x-9)^5}{27(x\!-\!1)(125x\!+\!3)^3},\
\frac{3125\,x^2(x-1)^3(5x+27)^5}{4(625x^3\!-\!2875x^2\!+\!675x\!-\!729)^3}.
\end{equation}

\begin{table}
\begin{center} \begin{tabular}{|c|c|c|c|c|c|c|c|c|}
\hline $\!$Denominators$\!$ & Schwarz & Klein &
\multicolumn{4}{|c|}{$g(D)$, when $\deg\gamma:\widetilde{D}\!\to D$ is equal to} \\
\cline{4-7} $k,\ell,m$ & type & degree & $\quad\;m\;\quad$
& $\quad\;\ell\;\quad$ & $\quad\;k\;\quad$ & $\quad\;1\;\quad$ \\
\hline 2, 3, 3 & $(1/2,\,1/3,\,1/3)$ & 1 & 0 & 0 & 0 & 0 \negvspace\\
& $(1/3,\,1/3,\,2/3)$ & 2& 0 & 0 & 1 & 1 \negvspace\\
2, 3, 4 & $(1/2,\,1/3,\,1/4)$ & 1& 0 & 0 & 0 & 0 \negvspace\\
& $(2/3,\,1/4,\,1/4)$ & 2& 0 & 1 & 2 & 3 \negvspace\\
2, 3, 5 & $(1/2,\,1/3,\,1/5)$ & 1 & 0 & 0 & 0 & 0 \negvspace\\
& $(1/2,\,1/3,\,2/5)$ & 7 & 0 & 0 & 0 & 0 \negvspace\\
& $(1/2,\,1/5,\,2/5)$ & 3 & 0 & 2 & 2 & 4 \negvspace\\
& $(1/3,\,1/3,\,2/5)$ & 2 & 1 & 1 & 3 & 5 \negvspace\\
& $(1/3,\,2/3,\,1/5)$ & 6 & 1 & 1 & 3 & 5 \negvspace\\
& $(2/3,\,1/5,\,1/5)$ & 2 & 1 & 3 & 5 & 9 \negvspace\\
& $(1/3,\,2/5,\,3/5)$ & 10 & 1 & 3 & 5 & 9 \negvspace\\
& $(1/3,\,1/5,\,3/5)$ & 4 & 1 & 3 & 5 & 9 \negvspace\\
& $(1/5,\,1/5,\,4/5)$ & 6 & 1 & 5 & 7 & 13 \negvspace\\
& $(2/5,\,2/5,\,2/5)$ & 6 & 1 & 5 & 7 & 13 \negvspace\\ \hline
\end{tabular} \end{center}
\caption{Genus of Darboux curves} \label{gentab}
\end{table}
Table \ref{gentab} shows the genus of all Darboux curves; this information
can be computed using Lemma \ref{dfiberprod}, part (iii) of Lemma
\ref{dramifico}, remarks in Appendixes \ref{fiberprod} and
\ref{apullbacks}. 
The third column shows the degree of Klein morphism for the main
representatives of each Schwarz type. We explicitly compute and utilize
Darboux coverings of minimal degree.

\subsection{Hypergeometric functions on Darboux curves}

The purpose of this subsection is to describe technical features of 
identifying and representing radical functions. 
First we describe radical solutions of a Fucshian equation with cyclic
monodromy.
\begin{lemma} \label{cctwosols}
Let $C$ be an algebraic curve. 
Consider a second order Fuchsian differential equation $(\ref{genlhde})$ on
$C$. Suppose that its differential Galois group $G$ is finite and cyclic.
Suppose that there are singular points with non-integer local exponent
differences. Then:
\begin{enumerate}
\item There exist exactly two independent radical solutions $f_1,f_2$.%
\item For any regular singular point $P\in C$ of $(\ref{genlhde})$ where the
exponent difference is not an integer, the 
local exponents of $f_1$ and $f_2$ are different.
\end{enumerate}
\end{lemma}
\proof Let $K_{PV}\supset\CC(C)$ be the Piccard-Vessiot extension of the
differential equation (\ref{genlhde}), and let $V$ denote its space of
solutions in $K_{PV}$. The monodromy group $G$ acts does not act on $V$ by
scalar multiplication, because otherwise the quotient of two independent
solutions would be in $\CC(C)$. Hence the representation of $G$ on $V$
splits into two irreducible representations. Let $f_1,f_2\in V$ be
generators of the two $G$-invariant subspaces. If $P\in C$ is a point where
the exponent difference is not an integer, $f_1,f_2$ are in different
spaces invariant under the the local monodromy group. Hence they have
well-definedlocal exponents, which should be different. \qed

Hypergeometric functions which pull-back to radical solutions under a
Darboux covering 
are characterized as follows.
\begin{lemma} \label{radpull}
Let $H$ denote a hypergeometric equation. Suppose that its monodromy 
$G$ is tetrahedral, octahedral or icosahedral. Let $\phi:D\to\PP^1$ denote a
Darboux morphism for $H$, of degree $d$. 
Let $\theta(z)\,\hpg21{A,B}{C}{f(z)}$ denote a hypergeometric solution
of $H$, with $f(z)$ a fractional-linear function, 
and $\theta(z)$ a radical factor. Assume that the denominator of the lower
parameter $C$ is equal to $|G|/d$. Then the hypergeometric function is
pull-backed to a radical function. Conversely, each radical solution of the
pull-backed eqaution represents (up to a scalar multiple) a hypergeometric
equation with the assumed property.
\end{lemma}
\proof Set $m=|G|/d$. Suppose that the pull-backed equation is normalized so
that its monodromy is the cyclic group of order $m$. Then all
non-integer local exponents of the pull-backed equation have denominator
$m$. We may assume that $f(z)=z$. There is a point $P$ above $z=0$ such that
the exponent difference $\delta$ for the pull-backed equation at $P$
is non-integer. The denominator of $\delta$ is $m$. Let
$\lambda_1,\lambda_2$ denote the local exponents at $P$, and let $t$ denote
a local parameter at $P$. For each local exponent $\lambda_j$ there is a
unique power series solution of the form
$t^{\lambda_j}(1+\alpha_1t+\alpha_2t+\ldots)$. By part {\em (ii)} of Lemma
\ref{cctwosols}, both power series represent radical functions.

The hypergeometric function is pull-backed, up to a constant multiple, to
one of the two mentioned power series. Hence the pull-back is a radical
function. Conversely, push-forwards of the two power series must me
hypergeometric series.\qed
%
\begin{lemma} \label{contigquo}
Let $F_1$, $F_2$ denote two contiguous algebraic (but not rational) Gauss
hypergeometric functions. Let $D$ denote their common Darboux curve. Then
the pull-back of $F_1/F_2$ is a rational function on $D$.
\end{lemma}
\proof The logarithmic derivative $F_1'/F_1$ is a rational function $D$. Up
to a factor in $C(D)$, $F_1'$ is Gauss hypergeometric function contiguous to
$F_1$. The contiguous relation between $F_1$, $F_1'$ and $F_2$ has
coefficients in $\CC(x)$, hence the claim follows. \qed

If we take any pair of contiguous evaluations in 
\S \ref{database}, the quotient of the right-hand sides is a rational (rather
then radical) function on the Darboux curve, as suggested by Lemma
\ref{contigquo}. This allows us to express other contiguous evaluations
conveniently as product of a fixed radical function and some rational
function on the Darboux curve.


\section{Computation of Darboux evaluations}
 \label{databasec}

This is an outline of actual computations that led us to the list of Darboux
evaluations in \S \ref{database}. If the Darboux curve has genus 0,
computations are quite straightforward. Subsections \ref{compfuncs} and
\ref{compevals} are devoted to difficulties of expressing rational and
radical functions on genus 1 Darboux curves.

\subsection{Computation of Darboux curves and coverings}
\label{darbcmorphs}

Darboux coverings for the three standard Schwarz types are considered in
\S \ref{darboux4st}. To compute other Darboux coverings and curves,
Lemma \ref{dfiberprod} is used . Let us fix an Schwarz type
which is not standard, and let $E_1$ denote its representative hypergeometric
equation as listed in \S \ref{alghypers}. Suppose that
$\varphi:\PP_z^1\mapsto\PP_Z^1$ is a Darboux covering (of minimal degree)
for the corresponding standard hypergeometric equation. Suppose that
$\psi:\PP_X^1\mapsto\PP_Z^1$ is a Klein covering for $E_1$. Then the Darboux
curve is a fiber product of $\varphi$ and $\psi$. An equation for it is
given by the equation $\psi(X)=\phi(z)$; see \S \ref{fiberprod}.

If the Darboux curve has genus 0, then a parameterization of it immediately
gives the Darboux covering. As an example, consider the icosahedral type
$(1/2,1/3,2/5)$. We apply Lemma \ref{dfiberprod} with $E_0=E(1/5,1/2,1/3)$
and $E_1=E(2/5,1/2,1/3)$. By formulas (\ref{isophi1}) and (\ref{okleinms})
we get the following equation for the Darboux curve:
\begin{equation}
\frac{X^2\,(189-64X)^5}{(3584X^2\!+\!2457X\!-\!2916)^3}=
\frac{1728\,z\,(z^2-11z-1)^5}{(z^4+228z^3+494z^2-228z+1)^3}.
\end{equation}
This is a rational curve. It can be parameterized by standard algorithms and
computer algebra systems such as {\sf Maple}. Here is a parameterization:
\begin{equation}
X=\frac{1728\,x\,(x^2-11x-1)^5}{(x^4+228x^3+494x^2-228x+1)^3},\qquad
z=\frac{(7x-1)^5}{x^2\,(x+7)^5}.
\end{equation}
(Recall that parameterizations are unique up to fractional-linear
transformations on $\PP_x^1$.) The parametric expression for $X$ gives the
Darboux covering for $E(2/5,1/2,1/3)$. We recognize that this is the same
Darboux covering (\ref{isophi1}) as for the Schwarz type $(1/2,1/3,1/5)$.

With a computer algebra package and standard algorithms \cite{mvhoeij},
\cite{kovacic} at hand, it is straightforward to pull-back a hypergeometric
equation to a Darboux curve of genus 0, and solve the pull-backed
differential equation. Lemma \ref{radpull} characterizes hypergeometric
functions which have to be identified with radical solutions. Since there
are only 2 radical solutions by Lemma \ref{cctwosols}, computer algebra
systems should return them. (Otherwise we may consider a simplified version
of the procedure in \S \ref{compevals}.)

If the Darboux curve $D$ has genus 1, we first wish to compute a convenient
Weierstrass model from the equation $\psi(X)=\varphi(z)$. With such a model
at hand, we identify $D$ with the elliptic curve $(D,\mathcal{O})$, where
$\mathcal{O}$ denotes the point at infinity. We always seek to have the point
$(x,\xi)=(0,0)$ on $D$ above $X=0$, so to allow easy power series
verification of the evaluations in (\ref{icosellipta})--(\ref{icoselliptz}).

The Darbox covering is given by the $X$-component of an isomorphism between
the elliptic curve and the model $\psi(X)=\varphi(z)$. We wish to express
the covering function in such a way that multiplicities of its zeroes and
poles would be well visible. This is not a straightforward problem; it is
discussed in \S \ref{compfuncs}. Pull-backs of hypergeometric
equations onto elliptic curves and finding radical solutions of those
pull-backs are discussed in \S \ref{compevals}. In the rest of this
subsection, we derive the four elliptic curves $C_3$, $C_4$, $C_5$, $C_6$,
introduced in (\ref{darbouxc1}), (\ref{darbouxc2}), (\ref{darbouxc3}),
(\ref{darbouxc4}) as Darboux curves for some icosahedral Schwarz types.

For the Schwarz type $(1/3,1/3,2/5)$, we use Lemma \ref{dfiberprod} with
$E_0=E(1/5,1/2,1/3)$ and $E_1=E(2/5,1/3,1/3)$. This gives the equation
$X^2\big/4(X\!-\!1)=\varphi_1(z)$ for the Darboux curve. After applying the
fractional-linear transformation $F\mapsto F/(F-1)$ to both sides we get:
\begin{equation} \label{fibprod3}
\frac{X^2}{(X-2)^2}=-\frac{1728\;z\;(z^2-11z-1)^5}
{(z^2+1)^2\,(z^4-522z^3-10006z^2+522z+1)^2}.
\end{equation}
We collect full squares onto the left-hand side and observe that the Darboux
curve is isomorphic to the genus 1 curve
\begin{equation}
\widetilde{\xi}{}^2=-1728\,z\,(z^2-11z-1). 
\end{equation}
This curve is isomorphic to $C_3$ via the isomorphism
$(z,\widetilde{\xi})\mapsto (3x,72\xi)$. The Darboux covering is given by
the $X$-component of the isomorphism between $C_3$ and (\ref{fibprod3}). We
have
\[
\frac{X}{X-2}=
\frac{72\,\xi\,(9x^2-33x-1)^2}{(9x^2\!+1)(81x^4-14094x^3-90054x^2+1566x+1)},
\]
so
\begin{equation} \label{darbexp3}
\varphi_3(x,\xi)= 2\Big/ 
\left(1+\frac{(9x^2\!+1)(81x^4-14094x^3-90054x^2+1566x+1)}
{72\,\xi\,(9x^2-33x-1)^2}\right). 
\end{equation}
Expression (\ref{isophi3}) is derived by methods of \S \ref{compfuncs}.

For the Schwarz type $(1/3,2/3,1/5)$, we use Lemma \ref{dfiberprod} with
$E_0=E(2/5,1/2,1/3)$ and $E_1=E(4/5,1/3,1/3)$. Then we get the same equation
$X^2\big/4(X\!-\!1)=\varphi_1(z)$. Hence the Darboux curve is the same as
for the type $(1/3,1/3,2/5)$.

For the Schwarz type $(2/3,1/5,1/5)$, we use Lemma \ref{dfiberprod} with
$E_0=E(1/5,1/2,1/3)$ and $E_1=E(1/5,1/5,2/3)$. This gives the equation
$4X(1-X)=\varphi_1(z)$ for the Darboux curve. After applying the
fractional-linear transformation $F\mapsto 1-F$ to both sides we get:
\begin{equation} \label{fibprod4}
(1-2X)^2=\frac
{(z^2+1)^2\,(z^4-522z^3-10006z^2+522z+1)^2}{(z^4+228z^3+494z^2-228z+1)^3}.
\end{equation}
We collect full squares to the left-hand side and observe that the Darboux
curve is isomorphic to the genus 1 curve
\begin{equation} \label{fibprod4a}
\widetilde{\xi}{}^2=z^4+228z^3+494z^2-228z+1.
\end{equation}
This curve is isomorphic to $C_4$ by the isomorphism
\begin{equation} \label{fibprod4b}
(z,\widetilde{\xi})\mapsto \left( \frac{57x-5\xi}{x+25},\,
\frac{25\,(\xi^2-570\xi-380x^2+248x+25)}{(x+25)^2}\right).
\end{equation}
Like with $\varphi_3(x,\xi)$, we 
identify $\varphi_4(x,\xi)$ with the $X$-component of the isomorphism
between $C_4$ and (\ref{fibprod4}), and apply methods of Lemma
\ref{compfuncs} to get expression (\ref{isophi4}).

For the Schwarz type $(1/3,2/5,3/5)$, we use Lemma \ref{dfiberprod} with
$E_0=E(2/5,1/2,1/3)$ and $E_1=E(2/5,2/5,2/3)$. Then we get the same equation
$4X(1-X)=\varphi_1(z)$. Hence the Darboux curve is the same as for the type
$(2/3,1/5,1/5)$.

For the Schwarz type $(1/3,1/5,3/5)$, we use Lemma \ref{dfiberprod} with
$E_0=E(1/5,1/2,1/3)$ and $E_1=E(1/5,1/3,3/5)$. This gives the equation
$-64X\big/(X-1)(9X-1)^3=\varphi_1(z)$ for the Darboux curve. This curve is
isomorphic to $C_5$, though it is not straightforward to compute a handy
isomorphism with current computer algebra packages. The package {\sf
algcurves} of {\sf Maple 9.0} can be used to obtain a Weierstrass form and
an isomorphism. The isomorphism ought to be simplified using methods of
\S \ref{compfuncs}. Eventually, we obtain an isomorphism given by
$z=-2(\xi-3x)^2\big/(\xi+3x)(4x+1)$ and $X=\varphi_5(x,\xi)$ as in
(\ref{isophi5}).

For the Schwarz type $(1/5,1/5,4/5)$, we use Lemma \ref{dfiberprod} with
$E_0=E(1/5,1/2,2/5)$ and $E_1=E(1/5,1/5,4/5)$. This gives the equation
$4X(1-X)=\varphi_2(z)$ for the Darboux curve. After applying the
fractional-linear transformation $F\mapsto 1-F$ to both sides we get
\begin{equation}
(1-2X)^2=-\frac{64\,x\,(x^2-x-1)^5}{(x^2+1)^2(x^4-22x^3-6x^2+22x+1)^2}.
\end{equation}
We collect full squares to the left-hand side and easily observe that the
Darboux curve is isomorphic $C_6$.

For the Schwarz type $(2/5,2/5,2/5)$, we use Lemma \ref{dfiberprod} with
$E_0=E(2/5,1/2,1/5)$ and $E_1=E(2/5,2/5,2/5)$. Then we get the same equation
$4X(1-X)=\varphi_2(z)$. Hence the Darboux curve is the same as for the type
$(1/5,1/5,4/5)$.

\begin{table}
\begin{center} \begin{tabular}{|c|c|c|}
\hline Elliptic & Mordell-Weil & Rational \vspace{-3pt}\\
curve & group &  points \\ \hline
 $C_3$ & $\ZZ/6\ZZ$  & $\mathcal{O}$,
$(0,0)$, $(-\frac19,\frac59)$,
$(-\frac19,-\frac59)$, $(1,-5)$, $(1,5)$\\
$C_4$ & $\ZZ\oplus\ZZ/2\ZZ$  & $\mathcal{O}$, $\mathcal{O}^*\!=(0,0)$,
$A_n=n\left(\frac15,\frac35\right)$,
$\widetilde{A}_n=-n\left(\frac15,\frac35\right)$,\\
&&$A^*_n=n\left(\frac15,\frac35\right)+\mathcal{O}^*$,
$\widetilde{A}^*_n=-n\left(\frac15,\frac35\right)+\mathcal{O}^*$\\
$C_5$ & $\ZZ/4\ZZ\oplus\ZZ/2\ZZ$  & $\mathcal{O}$, $(0,0)$,
$(-\frac14,-\frac34)$, $(-\frac14,-\frac34)$, \\ 
&& $(-1,0)$, $(-\frac1{16},0)$,
 $(\frac14,-\frac54)$, $(\frac14,\frac54)$\\
$C_6$& $\ZZ/6\ZZ$ & $\mathcal{O}$, $(0,0)$, $(-1,1)$, $(-1,-1)$,
$(1,-1)$, $(1,1)$\\
\hline
\end{tabular} \end{center}
\caption{Rational points on elliptic curves} \label{ratpoints}
\end{table}

Notice that the elliptic curves $C_3$, $C_4$, $C_5$, $C_6$ are defined over
$\QQ$. It is useful to know rational points on them. Table \ref{ratpoints}
gives this arithmetic information \cite{silverman1}. It was computed using
{\sf Maple} package {\sf Apecs} \cite{apecs}. Recall that by $\mathcal{O}$
we denote the point at infinity. As we see, only the curve $C_4$ has
infinitely many rational
points. In Table \ref{ratpoints} we introduce the notation 
$\mathcal{O}^*$, $A_n$, $\widetilde{A}_n$, $A_n^*$, $\widetilde{A}_n^*$
(with positive $n\in\ZZ$) for the rational points on $C_4$.

\subsection{Representing functions on genus 1 curves}
\label{compfuncs}

Here we consider the problem of representation of rational functions on
elliptic curves. Foremost, we use techniques of this subsection to compute
expressions (\ref{isophi3}), (\ref{isophi4}), (\ref{isophi5}),
(\ref{isophi6}) for Darboux coverings from genus 1 Darboux curves.
Subsection \ref{compevals} extends these techniques for computation of
expressions on the right-hand sides of
(\ref{icosellipta})--(\ref{icoselliptz3}),
(\ref{icosellipta4})--(\ref{icoselliptz4}),
(\ref{icosellipta5})--(\ref{icoselliptz5}),
(\ref{icosellipta6})--(\ref{icoselliptz}).

A canonical way to represent a rational functions $F$ on a (hyper)elliptic
curve \mbox{$\xi^2=G(x)$} (with $G(x)\in\CC[x]$) is the sum $f_1(x)+\xi
f_2(x)$, with $f_1(x),f_2(x)\in\CC(x)$. This representation suits well
algebraic computations, but it gives little geometric information about the
function. For example, the principal divisor for a function can be much
simpler than the degree of $f_1(x)$ and $f_2(x)$ may suggest. We would like
to have a compact expression that reflects well multiplicities in the
principal divisor. We do not give strict definitions or algorithms for an
alternative representation. Rather, we give tables of principal divisors on
the elliptic curves $C_3$, $C_4$, $C_5$, $C_6$, and propose to combine those
principal divisors to make the divisor for $F$. The corresponding
multiplicative expression in $\CC[x,\xi]$-polynomials from our tables will
give, up to a constant multiple, a compact expression for $F$ that we 
seek. We need to compute only finitely many rational and radical
functions on elliptic curves, and our tables give enough information
for these purposes.

Concretely, we start with Darboux covering $\varphi_3(x,\xi)$. Its principal
divisor can be computed from (\ref{darbexp3}) to be the following:
\begin{equation} \label{phi3div} \textstyle
(0,0)+\mathcal{O}+5\left( \frac{11+5\sqrt{5}}{6}, 0 \right) +5\left(
\frac{11-5\sqrt{5}}{6}, 0 \right)-3R_1-3R_2-3R_3-3R_4.
\end{equation}
Here the four points $R_1,R_2,R_3,R_4$ are defined by the equations
\begin{equation} \label{denomeq}
81x^4+6156x^3+4446x^2-684x+1=0, \qquad 150\xi=27x^3+1989x^2+741x-7.
\end{equation}
Table \ref{divisor1} gives a list of principal divisors on $C_3$.
\begin{table}
\[\begin{array}{|c|c|} \hline
\mbox{Function} & \mbox{Divisor} \\ \hline\vspace{-11pt} & \\
\xi & (0,0)+\left( \frac{11+5\sqrt{5}}{6}, 0 \right)+\left(
\frac{11-5\sqrt{5}}{6}, 0 \right)-3\mathcal{O} \\
1+33x-9x^2 & 2\left( \frac{11+5\sqrt{5}}{6}, 0 \right)+
2\left( \frac{11-5\sqrt{5}}{6}, 0 \right)-4\mathcal{O} \\
1-9\xi+54x & 3\left(-\frac19,-\frac59\right)-3\mathcal{O} \\
1+9\xi+54x & 3\left(-\frac19,\frac59\right)-3\mathcal{O} \\
1+9x & \left(-\frac19,-\frac59\right)+
\left(-\frac19,\frac59\right)-2\mathcal{O} \\
\xi+5x & (0,0)+(-\frac19,\frac59)+(1,-5)-\!3\mathcal{O}\\
1\!+\!21\xi\!-\!117x\!+\!9x\xi\!-\!234x^2 &
R_1+R_2+R_3+R_4+\left(-\frac19,-\frac59\right)-5\mathcal{O} \\
1\!-\!21\xi\!-\!117x\!-\!9x\xi\!-\!234x^2 &
\widetilde{R}_1+\widetilde{R}_2+\widetilde{R}_3+\widetilde{R}_4
+\left(-\frac19,\frac59\right)-5\mathcal{O} \\\hline\end{array}
\] \caption{Principal divisors on $C_3$} \label{divisor1}
\end{table}
For $i\in\{1,2,3,4\}$, by $\widetilde{R}_i$ we denote the inverse of $R_i$
in the group structure of $C_3$. Divisor (\ref{phi3div}) can be rewritten as
follows:
\begin{eqnarray*} \textstyle
\left\{ (0,0)+\!\left( \frac{11+5\sqrt{5}}{6}, 0 \right)\!+\!\left(
\frac{11-5\sqrt{5}}{6}, 0 \right) \!-\! 3\mathcal{O} \right\} + 2\left\{
2\left( \frac{11+5\sqrt{5}}{6}, 0 \right)\!+\!2\left(\frac{11-5\sqrt{5}}{6},
0\right)\!-\!4\mathcal{O}\right\}\nonumber \\
\textstyle +\Big\{ 3\left(-\frac19,-\frac59\right)-3\mathcal{O} \Big\}
-3\,\Big\{ R_1+R_2+R_3+R_4+
\left(-\frac19,-\frac59\right)-5\mathcal{O}\Big\}.
\end{eqnarray*}
Observe that each divisor in curled brackets is present in Table
\ref{divisor1}. We can immediately build the corresponding multiplicative
combination of the functions $\xi$, $1+33x-9x^2$, $1-9\xi+54x$,
$1+21\xi-117x+9x\xi-234x^2$. Up to undetermined constant multiple, the
multiplicative expression is (\ref{isophi3}). The constant multiple can be
determined by evaluating the multiplicative expression and (\ref{darbexp3})
a convenient point, say $(-\frac19,\frac59)$.

As an extra exercise, one may consider the function $1-\varphi_3$. Its
divisor can be computed from (\ref{isophi3}) or (\ref{darbexp3}) to be
\begin{equation} \label{phi3m1div}
3\widetilde{R}_1+3\widetilde{R}_2+3\widetilde{R}_3
+3\widetilde{R}_4-3R_1-3R_2-3R_3-3R_4.
\end{equation}
A straightforward combinatorial work suggests the expression
\begin{equation}
1-\varphi_3(x,\xi)=\frac{(1-21\xi-117x-9x\xi-234x^2)^3\,(1-9\xi+54x)}
{(1+21\xi-117x+9x\xi-234x^2)^3\,(1+9\xi+54x)}.
\end{equation}

Our proposal boils down in building a sufficient table of principal
divisors, and combining the known principal divisors to arrive at the
principal divisor of a target function. In practise, both things are done in
parallel. We start with the functions we wish to express, and compute their
divisors. We look at $\QQ$-rational points that occur, and use knowledge of
the Mordell-Weil group (see Table \ref{ratpoints}) to foresee and compute
suitable $\CC[x,\xi]$-polynomials that vanish only on rational points. Then
we distinguish $\QQ$-irreducible divisor components of higher degree. For
each such divisor component $\Gamma$, we use Gr\"{o}bner bases to find
$\CC[x,\xi]$-polynomials of minimal degree that vanish on $\Gamma$ with
sufficient multiplicities. We choose those polynomials whose divisors
enlarge $\Gamma$ minimally or least awkwardly. We look at additional
components that occur (usually they are rational points); if they are new,
we introduce new polynomials that could compensate the additional
components.

For example, consider the $\QQ$-irreducible component $R_1+R_2+R_3+R_4$
defined by (\ref{denomeq}). A Gr\"{o}bner basis gives the following
quadratic polynomials that vanish on it:
\begin{equation} \label{div33}
9x\xi-234x^2+21\xi-117x+1, \quad 3\xi^2-2088x^2+150\xi-744x+7.
\end{equation}
Other quadratic polynomials are obtained by linear combination. We have
chosen the first polynomial in (\ref{div33}), and consequently we had to
compensate its additional component $(-\frac19,-\frac59)$. A reasonable
alternative is the quadratic polynomial $\xi^2-21x\xi-150x^2+\xi+25x$, whose
divisor is $R_1+R_2+R_3+R_4+(0,0)+(1,-5)-6\mathcal{O}$.

Next consider computation of expression (\ref{isophi4}) for the Darboux
covering $\varphi_4(x,\xi)$. This is the most complicated case, so our
description of the computational method reaches deeper refinement level. Let
us introduce the functions $N$ and $L$ on the rational points of $C_4$:
\begin{equation} \label{defnl}
\begin{array}{ccccc}
N(A_n)=n, & N(A^*_n)=n, & N(\widetilde{A}_n)=-n,&
N(\widetilde{A}^*_n)=-n,& N(\mathcal{O}^*)=0,\\
L(A_n)=0,& L(A^*_n)=1, & L(\widetilde{A}_n)=0,& L(\widetilde{A}^*_n)=1,&
L(\mathcal{O}^*)=1.
\end{array}
\end{equation}
Principal divisors of $\CC[x,\xi]$-polynomials have the form $\sum_{j=1}^n
S_j-n\mathcal{O}$. If all points $S_j$ are rational, by Lemma \ref{princdiv}
we must have
\begin{equation} \label{pridivco}
\sum_{j=1}^n N(S_j)=0\quad\mbox{and}\quad \sum_{j=1}^n
L(S_j)\quad\mbox{even.}
\end{equation}

A preliminary expression for $\varphi_4(x,\xi)$ can be computed by composing
the obvious isomorphism between the curves in (\ref{fibprod4}) and
(\ref{fibprod4a}) with isomorphism (\ref{fibprod4b}). The principal divisor
of $\varphi_4(x,\xi)$ is:
\begin{equation} \label{phi4div} \textstyle
(0,0)+\mathcal{O}+5P_1+5P_2 
-\left(\frac1{125},\frac{57}{625}\right)-(-25,285)-5Q_1-5Q_2,
\end{equation}
where
\[ \textstyle
P_{1,2}=\left( -\frac32\pm\frac7{2\sqrt{5}}, 7\mp3\sqrt{5} \right), \qquad
Q_{1,2}=\left( -4\pm\frac9{\sqrt{5}}, -27\pm12\sqrt{5} \right).
\]
In the notation of Table \ref{ratpoints}, we have
$\left(\frac1{125},\frac{57}{625}\right)=A_5$ and $(-25,285)=A_5^*$.

\begin{table}
\[
\begin{array}{|c|c||c|c|} \hline
\mbox{Function} & \mbox{Divisor} & \mbox{Function} &
\mbox{Divisor} \\ \hline\vspace{-14pt} &&& \\
x & 2\,\mathcal{O}^*-2\mathcal{O} &
1-15x-5x^2 & P_1+P_2+\widetilde{P}_1+\widetilde{P}_2-4\mathcal{O} \\
1-5x & A_1+\widetilde{A}_1-2\mathcal{O} &
1-40x-5x^2 & Q_1\!+\!Q_2\!+\!\widetilde{Q}_1\!+\!\widetilde{Q}_2\!-\!4\mathcal{O} \\
25+x & A_5^*+\widetilde{A}_5^*-2\mathcal{O} & 5\!-\!7\xi\!-\!45x\!-\!5x^2 &
P_1+P_2+\widetilde{A}^*_5+
\widetilde{A}_1-4\mathcal{O} \\
1-125x & A_5+\widetilde{A}_5-2\mathcal{O} & 5\!+\!18\xi\!-\!80x\!+\!5x^2 &
Q_1+Q_2+\widetilde{A}^*_5+
A_1-4\mathcal{O} \\
5\xi+57x & A^*_5\!+\!\widetilde{A}_5\!+\!\mathcal{O}^*-3\mathcal{O} &
1-7y+15x+15x^2 &
P_1+P_2+\widetilde{A}_4^*+\widetilde{A}_2-4\mathcal{O} \\
-5\xi+57x & \widetilde{A}^*_5\!+\!A_5\!+\!\mathcal{O}^*-3\mathcal{O} &
4-7\xi-30x & P_1\!+\!P_2+\!\widetilde{A}^*_6-3\mathcal{O}\\
1\!+\!5\xi\!+\!10x & 2\widetilde{A}_1\!+\!A_2-3\mathcal{O} &
1+3\xi-20x & Q_1\!+\!Q_2\!+\!\widetilde{A}^*_4-3\mathcal{O} \\
1-3\xi+2x & A_2\!+\!A^*_2\!+\!\widetilde{A}^*_4-3\mathcal{O} &
1\!-\!8\xi\!+\!22x\!-\!15x^2 &
2A_1+A_2^*+\widetilde{A}^*_4-4\mathcal{O} \\
4\!+\!21\xi\!+\!41x & A_2\!+\!A^*_4\!+\!\widetilde{A}^*_6-3\mathcal{O} &
4-35\xi-101x &
A_1^*\!+\!A_5\!+\!\widetilde{A}^*_6-3\mathcal{O} \\
5-3\xi-34x & A^*_5\!+\!\widetilde{A}^*_4\!+\!\widetilde{A}_1-3\mathcal{O} &
20-7\xi-79x  & A_1\!+\!A^*_5\!+\!\widetilde{A}^*_6-3\mathcal{O}
\\  \hline
\end{array}\vspace{-3pt}
\]
\[
\begin{array}{|c|c|}
\hline \mbox{Function} & \mbox{Divisor} \\ \hline\vspace{-14pt} & \\
1\!+\!50x\!-\!125\xi^2\!+\!450x\xi\!-\!500x^2 & 5A_1+\widetilde{A}_5-6\mathcal{O} \\
1\!+\!50x\!-\!125\xi^2\!-\!450x\xi\!-\!500x^2 & 5\widetilde{A}_1+A_5-6\mathcal{O} \\
25\!-\!570\xi\!+\!248x\!+\!\xi^2\!-\!380x^2 &
R_1+R_2+R_3+R_4+2\widetilde{A}^*_5-6\mathcal{O} \\
4\!+\!95\xi\!+\!83x\!+\!21\xi^2\!-\!475x\xi\!+\!40x^2 &
R_1+R_2+R_3+R_4+\widetilde{A}^*_6+\widetilde{A}^*_4-6\mathcal{O} \\
 \hline
\end{array}\]\caption{Principal divisors on $C_4$} \label{divisor2}
\end{table}

\begin{table} \[
\begin{array}{|c|c|}
\hline \mbox{Function} & \mbox{Divisor} \\ \hline\vspace{-11pt} & \\
\xi+5x & 2\left(\frac14,-\frac54\right)+(0,0)-3\mathcal{O} \\
\xi-5x & 2\left(\frac14,\frac54\right)+(0,0)-3\mathcal{O} \\
\xi+3x & 2\left(-\frac14,\frac34\right)+(0,0)-3\mathcal{O} \\
\xi-3x & 2\left(-\frac14,-\frac34\right)+(0,0)-3\mathcal{O} \\
1+\xi+x & \left(-\frac14,-\frac34\right)+\left(\frac14,-\frac54\right)+(-1,0) \\
1-2\xi+6x & P_1+P_2+\left(\frac14,\frac54\right)-3\mathcal{O} \\
1+12x+16x^2 & P_1+P_2+\widetilde{P}_1\!+\!\widetilde{P}_2-4\mathcal{O} \\
1-2\xi-14x & Q_1+Q_2+\left(\frac14,-\frac54\right)-3\mathcal{O} \\
1-28x+16x^2 & Q_1+Q_2+\widetilde{Q}_1\!+\!\widetilde{Q}_2-4\mathcal{O} \\
1+8\xi-28x+8x\xi-104x^2 & R_1+R_2+R_3+R_4+\left(\frac14,\frac54\right)-5\mathcal{O} \\
\hline
\end{array}
\] \caption{Divisors on $C_5$} \label{divisor3}
\end{table}

It seems convenient to consider the lines through $P_1$, $P_2$ and through
$Q_1$, $Q_2$. Their equations are $7\xi=4-30x$ and $3\xi=20x-1$,
respectively. The third points on these two lines are, respectively,
$\widetilde{A}_6^*$ and $\widetilde{A}_4^*$. But after adding or subtracting
extra divisor terms $5\widetilde{A}^*_6$, $5\widetilde{A}^*_4$ in
(\ref{phi4div}), it is very cumbersome to compensate them due to
(\ref{pridivco}). The $\CC[x,\xi]$ polynomials for compensating principal
divisors are expected to have very large coefficients.

Rather than introducing the sub-expression $P_1+P_2+\widetilde{A}^*_6-3\mathcal{O}$ 
in (\ref{phi4div}), we may try to work with the principal divisor
$P_1+P_2+\widetilde{A}^*_4+\widetilde{A}_2-4\mathcal{O}$. Then
$\widetilde{A}^*_4$ is compensated automatically. Keeping in mind that
the divisors $P_1+P_2$ and $Q_1+Q_2$ are equivalent (respectively) 
to $A^*_6$ and $A^*_4$ for the purposes of restrictions (\ref{pridivco}),
we work out the
following expressions of (\ref{phi4div}) as sums of principal divisors:
\begin{eqnarray*}
5\left(P_1+P_2+\widetilde{A}^*_4+\widetilde{A}_2-4\mathcal{O}\right)
+\left(A_5+\widetilde{A}^*_5+\mathcal{O}^*-3\mathcal{O} \right)\\
-5\left( Q_1+Q_2+\widetilde{A}^*_4-3\mathcal{O} \right) -
\left(5\widetilde{A}_2+2A_5-7\mathcal{O}\right) - \left(
A^*_5+\widetilde{A}^*_5-2\mathcal{O} \right),\!\!
\end{eqnarray*}
and
\begin{eqnarray*}
5\left(P_1+P_2+\widetilde{A}^*_4+\widetilde{A}_2-4\mathcal{O}\right)
+\left( A_{10}+\widetilde{A}^*_{10}+\mathcal{O}^*-3\mathcal{O}\right)\\
-5\left( Q_1+Q_2+\widetilde{A}^*_4-3\mathcal{O} \right) -
\left(5\widetilde{A}_2+A_{10}-6\mathcal{O}\right) - \left(
A_5+A^*_5+\widetilde{A}^*_{10}-3\mathcal{O} \right).\!\!
\end{eqnarray*}
Now we can build a table of $\CC[x,\xi]$-polynomials of the involved
principal divisors. Eventually, the two decompositions of (\ref{phi4div})
give the following expressions for $\varphi_4(x,\xi)$:
\[
\frac{8208\,(1-7\xi+15x+15x^2)^5\;(-5\xi+57x)}{(1\!+\!3\xi\!-\!20x)^5
(59375x^2\xi\!+\!12350x\xi\!+\!475\xi\!-\!166250x^3\!-\!49875x^2\!-\!3800x\!-\!19)
(x\!+\!25)},
\]
and, respectively,
\[
\frac{432\,(1-7\xi+15x+15x^2)^5\,(49495\xi+292441x)}{(1\!+\!3\xi\!-\!20x)^5
(2375\xi^2\!+\!7400x\xi\!+\!1150\xi\!-\!42000x^2\!-\!10200x\!-\!521)
(2605\xi\!+\!29678x\!-\!475)}.
\]
Compensation of $5\widetilde{A}_2$ look quite awkward in both formulas.

As the last attempt, we introduce
$P_1+P_2+\widetilde{A}^*_5+\widetilde{A}_1-4\mathcal{O}$ and
$Q_1+Q_2+\widetilde{A}^*_5+A_1-4\mathcal{O}$ in (\ref{phi4div}), forgetting
divisors of the linear polynomials $7\xi+30x-4$ and $3\xi-20x+1$. A natural
effort to compensate $5\widetilde{A}_1$ and $5A_1$ leads to the following
expression decomposition of (\ref{phi4div}):
\begin{eqnarray*}
5\left(P_1+P_2+\widetilde{A}^*_5+\widetilde{A}_1-4\mathcal{O}\right)
+\left(5A_1+\widetilde{A}_5-6\mathcal{O}\right)+
\left(2\mathcal{O}^*-2\mathcal{O}\right)\\
-5\left(Q_1+Q_2+\widetilde{A}^*_5+A_1-4\mathcal{O} \right) -
\left(5\widetilde{A}_1+A_5-6\mathcal{O}\right) - \left(
A^*_5+\widetilde{A}_5+\mathcal{O}^*-3\mathcal{O} \right).\!\!
\end{eqnarray*}
This expression gives formula (\ref{isophi4}). The most convenient principal
divisors are listed in Table \ref{divisor2}. The points $R_1,R_2,R_3,R_4$ on
$C_4$ are the points in the fiber $z=1$ of $\varphi_4(x,\xi)$.

The principal divisor for $\varphi_5(x,\xi)$ on $C_5$ turns out to be
\begin{equation} \label{phi5div} \textstyle
(0,0)+\mathcal{O}+5P_1+5P_2-\left(-\frac14,\frac34\right)
-\left(-\frac14,-\frac34\right)-5Q_1-5Q_2,
\end{equation}
where
\[ \textstyle
P_{1,2}=\left( \frac{-3\pm\sqrt{5}}8, \frac{-5\pm3\sqrt{5}}8 \right), \qquad
Q_{1,2}=\left( \frac{7\pm3\sqrt{5}}8, \frac{-45\mp21\sqrt5}8 \right).
\]
As mentioned in \S \ref{darbcmorphs}, computation of $C_5$ and a preliminary
expression for $\varphi_5(x,\xi)$ was not straightforward. A natural effort
leads to the following expression decomposition of (\ref{phi5div})
into principal ideals:
\begin{eqnarray} \label{phi3divm} \textstyle
3\,\Big( 2\left(\frac14,-\frac54\right)+(0,0)-3\mathcal{O}\Big)
+5\,\Big(P_1+P_2+\left(\frac14,\frac54\right)-3\mathcal{O}\Big) \nonumber \\
\textstyle -\Big( \left(\frac14,\frac54\right)
+\left(\frac14,-\frac54\right)+\left(-\frac14,\frac34\right)
+\left(-\frac14,-\frac34\right)-4\mathcal{O}\Big)\nonumber \\
\textstyle -2\,\Big( 2\left(\frac14,\frac54\right)+(0,0)-3\mathcal{O}\Big)
-5\,\Big(Q_1+Q_2+\left(\frac14,-\frac54\right)-3\mathcal{O}\Big).\!\!
\end{eqnarray}
This expression gives formula (\ref{isophi5}). The most convenient principal
divisors are listed in Table \ref{divisor3}. The points $R_1,R_2,R_3,R_4$ on
$C_5$ are the points in the fiber $z=1$ of $\varphi_5(x,\xi)$.

\begin{table}
\[\begin{array}{|c|c||c|c|} \hline
\mbox{Function} & \mbox{Divisor} & \mbox{Function} &
\mbox{Divisor} \\ \hline\vspace{-11pt} &&& \\
\xi & (0,0)+P_1+P_2-3\mathcal{O} &
1+x-x^2 & 2P_1+2P_2-4\mathcal{O} \\
1-\xi & (-1,1)+2(1,1)-3\mathcal{O} &
1-4x-x^2 & Q_1\!+\!Q_2\!+\!R_1\!+\!R_2\!-\!4\mathcal{O} \\
1+\xi & (-1,-1)+2(1,-1)-3\mathcal{O} & \xi+2x+x^2 &
(0,0)+Q_1+Q_2-3\mathcal{O} \\
1+\xi+2x & 3(-1,1)-3\mathcal{O} & 1+\xi-2x &
Q_1+Q_2+(1,1)-3\mathcal{O} \\
1-\xi+2x & 3(-1,-1)-3\mathcal{O} & 1-\xi-2x &
R_1+R_2+(1,-1)-3\mathcal{O} \\
\hline\end{array}
\] \caption{Divisors on $C_6$} \label{divisor4}
\end{table}

The principal divisor for $\varphi_6(x,\xi)$ on $C_5$ turns out to be
\begin{equation} \label{phi6div} \textstyle
(0,0)+\mathcal{O}+5P_1+5P_2-\left(1,1\right) -\left(-1,1\right)-5Q_1-5Q_2.
\end{equation}
where
\[ \textstyle
P_{1,2}=\left( \frac{1\pm\sqrt{5}}2, 0\right), \qquad Q_{1,2}=\left(
-2\pm\sqrt5, -5\pm2\sqrt5 \right).
\]
We naturally arrive at the following expression for (\ref{phi6div}):
\begin{eqnarray*} \label{phi6divm} \textstyle
\Big(P_1+P_2+(0,0)-\mathcal{O}\Big) +2\,\Big(2P_1+2P_2-4\mathcal{O}\Big)
+2\,\Big(2(1,1)+(1,-1)-2\mathcal{O}\Big)\nonumber\\
-\Big(3(-1,1)-3\mathcal{O}\Big) -5\,\Big(Q_1+Q_2+(1,1)-3\mathcal{O}\Big).
\end{eqnarray*}
This expression gives formula (\ref{isophi6}). The most convenient principal
divisors are listed in Table \ref{divisor4}. The points $R_1,R_2$ on $C_6$
are in the fiber $z=1$ of $\varphi_6(x,\xi)$.

\subsection{Computation of hypergeometric evaluations}
\label{compevals}

In principle, evaluations (\ref{icosellipta})--(\ref{icoselliptz3}),
(\ref{icosellipta4})--(\ref{icoselliptz4}),
(\ref{icosellipta5})--(\ref{icoselliptz5}),
(\ref{icosellipta6})--(\ref{icoselliptz}) are computed by pulling-back
their hypergeometric equations onto $C_3$, $C_4$, $C_5$, $C_6$,
respectively, 
and finding radical solutions of the pull-backed differential equations. 
Their divisors will have coefficients in $\QQ$ rather than $\ZZ$, 
as described in Appendix \ref{secalgcurves}.
But standard computer algebra systems do not
handle differential equations on higher genus curves.

The pull-backed equations are cumbersome as we will see. On the other hand,
their singular points and exponent differences are easy to see,
as explained in Appendix \ref{fuchsian}. The local exponents tell us possible
coefficients in the principal divisors of the radical solutions. Possible
principal divisors are restricted by Lemma \ref{principald} and Lemma
\ref{cctwosols}. For simplest hypergeometric equations, we may end up with
just 2 possible principal divisors for radical solutions. Then we can find
and check those solutions without computing the pull-back equation explicitly.
Additional contiguous evaluations can be obtained by differentiating known
solutions and contiguous relations, while respecting Lemma \ref{contigquo}
and avoiding explicit computation with pull-back equation again.

Generally, we may have several candidates for radical solutions, which we
must check by substituting into the pull-back equation. If we have one
undetermined simple zero for a radical solution, its location can be
restricted by the arithmetical argument that the principal divisor should be
invariant under the Galois action of $\overline{\QQ}$.

We start with computation of evaluations
(\ref{icosellipta})--(\ref{icoselliptc3}). We use Riemann notation
(\ref{riemanp}) and consider the icosahedral hypergeometric equation with
the solution space
\[
P\left\{\begin{array}{ccc} 0 & 1 & \infty \\
0 & 0 & -1/30 \\ 2/5  & 1/3 & 3/10
\end{array} \; z \right\}.
\]
Its pull-back $z\mapsto\varphi_3(x,\xi)$, $y(z)\mapsto Y(\varphi_3(x,\xi))$
onto $C_3$ has the following regular singular points and local exponents:
\[
\left\{\begin{array}{cccccccc} (0,0) & \mathcal{O} & \left(
\frac{11+5\sqrt{5}}{6}, 0 \right) & \left( \frac{11-5\sqrt{5}}{6}, 0 \right)
& R_1 & R_2 & R_3 & R_4 \\
0 & 0 & 0 & 0 & -\frac1{10} & -\frac1{10} &
-\frac1{10} & -\frac1{10} \\
\frac25 & \frac25 & 2 & 2 & \frac9{10} & \frac9{10} & \frac9{10} &
\frac9{10}
\end{array} \right\}
\]
Condition {\em (i)} of Lemma \ref{principald} leaves only two candidates for
the divisors of radical solutions:
\begin{eqnarray} \textstyle \label{icodivsa}
\frac{2}{5}\,\mathcal{O}-\frac{1}{10}R_1-\frac{1}{10}R_2
-\frac{1}{10}R_3-\frac{1}{10}R_4, \qquad
\frac{2}{5}\,(0,0)-\frac{1}{10}R_1-\frac{1}{10}R_2
-\frac{1}{10}R_3-\frac{1}{10}R_4.
\end{eqnarray}
They must represent divisors of the two solutions of Lemma (\ref{cctwosols})
and, by identification of local solutions, the hypergeometric series
(\ref{icosellipta}) and (\ref{icoselliptc3}). Since the point $x=0$
corresponds to $(0,0)\in C_3$, the first divisor correspond to the series in
(\ref{icosellipta}). The right-hand side of that formula is easy to
construct from the divisor by using Table \ref{divisor1} and the known value
of left-hand side at $x=0$. Similarly, the second divisor in
(\ref{icodivsa}) implies the identity
\[
\varphi_3(x,\xi)^{2/5}\hpg{2}{1}{7/10,11/30}{7/5}{\varphi_3(x,\xi)}=
\frac{x^{1/5}\,(1-9\xi+54x)^{1/30}}{(1+21\xi-117x+9x\xi-234x^2)^{1/10}}.
\]
This gives formula (\ref{icoselliptc3}). We derived these identities without
even computing the pull-back differential equation on $C_3$.  

Next we consider computation of evaluations
(\ref{icoselliptk3})--(\ref{icoselliptz3}). The pull-back
$z\mapsto\varphi_3(x,\xi)$, $y(z)\mapsto Y(\varphi_3(x,\xi))$ of the
icosahedral hypergeometric equation with the solution space
\[
P\left\{\begin{array}{ccc} 0 & 1 & \infty \\
0 & 0 & -1/10 \\ 1/5  & 1/3 & 17/30
\end{array} \; z \right\}.
\]
is the following differential equation (with coefficients in a convenient
form):
\begin{eqnarray} \label{pullbackeq3}
Y''&\!\!\!+\!\!\!& \nonumber
\frac{3(3\!+\!47\xi\!+\!1974x\!-\!2051\xi^2\!+\!2676x\xi\!+\!54348x^2
\!+\!33\xi^3\!-\!1002x\xi^2\!-\!1548x^2\xi)}{10\,\xi^2\,(1+21\xi-117x+9x\xi-234x^2)}\,Y'\\
&\!\!\!+\!\!\!&\frac{51\,\xi\,(1-9\xi+54x)\,(1-21\xi-117x-9x\xi-234x^2)}
{25\,x^2\,(1+9x)\,(1+21\xi-117x+9x\xi-234x^2)^2}\,Y \ = \ 0.
\end{eqnarray}
Like in the previous case above, we know singularities and local exponents
of this equation without cumbersome computations. Here they are:
\[
\left\{\begin{array}{cccccc} (0,0) & \mathcal{O} & R_1 & R_2 & R_3& R_4 \\ 0
& 0 & -\frac3{10} & -\frac3{10} &
-\frac3{10} & -\frac3{10} \\
\frac15 & \frac15 & \frac{17}{10} & \frac{17}{10} & \frac{17}{10} &
\frac{17}{10}
\end{array} \right\}
\]
Condition {\em (i)} of Lemma \ref{principald} gives the following candidates
for the divisors of radical solutions:
\begin{eqnarray*} \textstyle \label{icodivsb}
\frac{1}{5}\,\mathcal{O}+X-\frac{3}{10}R_1-\frac{3}{10}R_2
-\frac{3}{10}R_3-\frac{3}{10}R_4, \quad
\frac{1}{5}\,(0,0)+Y-\frac{3}{10}R_1-\frac{3}{10}R_2
-\frac{3}{10}R_3-\frac{3}{10}R_4.
\end{eqnarray*}
Here $X$ and $Y$ are some regular points of (\ref{pullbackeq3}). By condition
{\em (ii)} and the additional statement of the same lemma, these should be
torsion points on $C_3$ defined over $\QQ$. The possibilities for $X$ and
$Y$ are: $(-\frac19,-\frac59)$, $(-\frac19,\frac59)$, $(1,-5)$, $(1,5)$.
This gives 8 possible divisors of a radical solution. For each possibility,
one has to construct a radical function (in any form) with that divisor, and
to check whether it is a solution of (\ref{pullbackeq3}). Alternatively,
candidate solutions can be expanded in power series around $(x,\xi)=(0,0)$
and compared with the hypergeometric series in (\ref{icoselliptk3}) and
(\ref{icoselliptm3}). Then one does not have to know explicit equation
(\ref{pullbackeq3}), but has to find enough power series terms of {\em all}
candidate solutions (so that the right candidates could be selected). It
turns out that actual solutions have $X=(-\frac19,-\frac59)$ and $Y=(1,-5)$.
These two solutions and evaluations (\ref{icoselliptk3}),
(\ref{icoselliptm3}) are expressed in a convenient form by using methods
of \S \ref{compfuncs}. 

Next we consider evaluation of formulas
(\ref{icosellipta4})--(\ref{icoselliptc4}). A pull-back of $E(1/5,1/5,2/3)$ 
has the following singularities and local exponents:
\[
\left\{\begin{array}{cccccccccc} \mathcal{O}^* & \mathcal{O} &
R_1 & R_2 & R_3 & R_4 & A_5 & A^*_5 & Q_1 & Q_2 \\
0 & 0 & 0 & 0 & 0 & 0 & -\frac1{30} & -\frac1{30} &
-\frac16 & -\frac16 \\
\frac15 & \frac15 & 2 & 2 & 2 & 2 & \frac16 & \frac16 & \frac56 & \frac56
\end{array} \right\}
\]
Condition {\em (i)} of Lemma \ref{principald} gives the following candidates
for the divisors on $C_4$ of radical solutions:
\[
\begin{array}{rr} 
-\frac{1}{6}Q_1-\frac{1}{6}Q_2+\frac{1}{6}A_5+\frac{1}{6}A^*_5,& 
\frac{1}{5}\,\mathcal{O}+\frac{1}{5}\,(0,0)-\frac{1}{6}Q_1-\frac{1}{6}Q_2
-\frac{1}{30}A_5-\frac{1}{30}A^*_5,\\ 
\frac{1}{5}\mathcal{O}-\frac{1}{6}Q_1-\frac{1}{6}Q_2+\frac{1}{6}A_5-\frac{1}{30}A^*_5,&
\frac{1}{5}(0,0)-\frac{1}{6}Q_1-\frac{1}{6}Q_2-\frac{1}{30}A_5+\frac{1}{6}A^*_5,\\
\frac{1}{5}\mathcal{O}-\frac{1}{6}Q_1-\frac{1}{6}Q_2-\frac{1}{30}A_5+\frac{1}{6}A^*_5,&
\frac{1}{5}(0,0)-\frac{1}{6}Q_1-\frac{1}{6}Q_2+\frac{1}{6}A_5-\frac{1}{30}A^*_5.
\end{array}
\]
The candidate divisors are grouped into possible pairs for a basis
of actual solutions, following Lemma \ref{cctwosols}. The divisors in the first pair
do not satisfy condition {\em (ii)} of Lemma \ref{principald}. To decide the
right pair, one may take one divisor from each of the one two pairs,
construct an expression for a corresponding radical function, and compare
its power series around $(0,0)$. It turns out, the last pair is the right one. Once we
have the right divisors for (\ref{icosellipta4}) and (\ref{icoselliptc4}),
we can proceed similarly as in the pervious cases. Application of the
methods of \S \ref{compfuncs}
may require some combinatorial creativeness. 
For example, here is a convenient splitting of the right divisor for
(\ref{icosellipta4}):
\[ \textstyle
\frac{1}{6}\left(\widetilde{A}^*_4\!+\!A^*_5\!+\!\widetilde{A}_1-3\mathcal{O}\right)
-\frac{1}{6}\left(Q_1\!+\!Q_2\!+\!\widetilde{A}^*_4-3\mathcal{O}\right)
-\frac{1}{30}\left(5\widetilde{A}_1\!+\!A_5-6\mathcal{O}\right).
\]

Now we consider evaluation of formulas
(\ref{icoselliptk4})--(\ref{icoselliptz4}). We comment only the most
complicated step of choosing the divisors of the actual solutions of the
pull-back of $E(2/5,3/5,1/3)$. The singularities and local exponents are the following:
\[
\left\{\begin{array}{cccccccc} \mathcal{O}^* & \mathcal{O} &
P_1 & P_2 & A_5 & A^*_5 & Q_1 & Q_2 \\
0 & 0 & 0 & 0 & -\frac16 & -\frac16 & -\frac56 & -\frac56 \\
\frac25 & \frac25 & 2 & 2 & \frac{13}{30} & \frac{13}{30} & \frac{13}6 &
\frac{13}6 \end{array} \right\}
\]
Condition {\em (i)} of Lemma \ref{principald} gives the following candidates
for the divisors of radical solutions:
\begin{eqnarray*}\textstyle
-\frac16\,A_5-\frac16\,A^*_5-\frac56\,Q_1-\frac56\,Q_2+2P_1,\hspace{-3.5cm}
&& \hspace{-2.2cm}\textstyle
-\frac16\,A_5-\frac16\,A^*_5-\frac56\,Q_1-\frac56\,Q_2+2P_2, \\
& -\frac16\,A_5-\frac16\,A^*_5-\frac56\,Q_1-\frac56\,Q_2+X+Y, &
X+Y\in\{A_5,A^*_5\} \\
&\frac25\mathcal{O}+\frac25\mathcal{O}^*+\frac{13}{30}\,A_5
+\frac{13}{30}\,A^*_5-\frac56\,Q_1-\frac56\,Q_2, \\
& \frac25\mathcal{O}+\frac{13}{30}\,A_5-\frac16\,A^*_5
-\frac56\,Q_1-\frac56\,Q_2+X, & X\in\{A_2,A^*_2\} \\
&\frac25\mathcal{O}-\frac16\,A_5+\frac{13}{30}\,A^*_5
-\frac56\,Q_1-\frac56\,Q_2+X, & X\in\{A_2,A^*_2\},\\
& \frac25\mathcal{O}^*+\frac{13}{30}\,A_5-\frac16\,A^*_5
-\frac56\,Q_1-\frac56\,Q_2+X, & X\in\{A_2,A^*_2\} \\
&\frac25\mathcal{O}^*-\frac16\,A_5+\frac{13}{30}\,A^*_5
-\frac56\,Q_1-\frac56\,Q_2+X, & X\in\{A_2,A^*_2\}.
\end{eqnarray*}
Here the restrictions on the additional points $X$, $Y$ follow from
condition {\em (ii)} and the additional statement of Lemma \ref{principald}.
The first three possibilities can be paired only with the fourth divisor as
(divisors of) functions $f_1$, $f_2$ of Lemma \ref{cctwosols}. To refute
them, one has to check only that a function with the fourth divisor (as the
principal divisor) is not a solution of the pull-back equation. Other
possibilities have to be paired and checked as we did for equation
(\ref{pullbackeq3}). The right divisors are these:
\[ \textstyle
\frac25\mathcal{O}+\frac{13}{30}\,A_5-\frac16\,A^*_5
-\frac56\,Q_1-\frac56\,Q_2+A^*_2, \quad
\frac25\mathcal{O}^*-\frac16\,A_5+\frac{13}{30}\,A^*_5
-\frac56\,Q_1-\frac56\,Q_2+A_2.
\]
Now we can compute (\ref{icoselliptk4}), (\ref{icoselliptm4}) like in the
previous cases, etc.

Computation of evaluations (\ref{icosellipta5})--(\ref{icoselliptz5}),
(\ref{icosellipta6})--(\ref{icoselliptz}) is similar and not more
complicated. First we find singularities and local exponents of the pull-backed
equation (of a corresponding main hypergeometric equation); make a list of
possible divisors for radical solutions of the pull-back; use Lemma
\ref{cctwosols} to make a short divisor list for necessary check, find
radical functions for the candidate divisors from the short list; compare
their power series around $(0,0)$ with the expansions of the hypergeometric
series; take the right divisors and find a convenient expression for their
functions (using methods of \S \ref{compfuncs}).

Differentiation and computation of contiguous evaluations to (\ref{icosellipta})--(\ref{icoselliptz})
quickly leads to large expressions. For example, an expression contiguous to 
 (\ref{icosellipta})  is
\begin{eqnarray}
\hpg{2}{1}{3/10,\,29/30}{3/5}{\varphi_3(x,\xi)}  \;\equal  \nonumber \\
& & \hspace{-2.1cm} \label{icoselliptb3}
\frac{(1\!+\!21\xi\!-\!117x\!+\!9x\xi\!-\!234x^2)^{9/10}(1\!+\!9x)^2(1\!+\!198x\!-\!99x^2)}
{(1-9\xi+54x)^{29/30}\,(1\!-\!21\xi\!-\!117x\!-\!9x\xi\!-\!234x^2)^2}.
\end{eqnarray}
The divisor for this function on $C_3$ is
\[ \textstyle
\frac9{10}R_1+\frac9{10}R_2+\frac9{10}R_3+\frac9{10}R_4
+S_1+S_2+\widetilde{S}_1+\widetilde{S}_2+\frac25\mathcal{O}
-2\widetilde{R}_1-2\widetilde{R}_2-2\widetilde{R}_3-2\widetilde{R}_4,
\]
where $S_1,S_2,\widetilde{S}_1,\widetilde{S}_2$ are the points with the
$x$-coordinate equal to $1\pm 10/3\sqrt{11}$.
For least painful computation of contiguous expressions, keep in mind that 
the quotient of contiguous hypergeometric functions is a rational function on a Darboux curve,
by Lemma \ref{contigquo}. In particular, the above expression should be viewed as 
$G(x,\xi)\times\hpg{2}{1}{3/10,\,-1/30}{3/5}{\varphi_3}$ with
$G(x,\xi)\in\CC(C_3)$. To compute $G(x,\xi)$,
only divisors with integer (rather than $\QQ$) coefficients have to be considered
as in \S \ref{compfuncs}.

\section{Appendix}
\label{appendix}

Here we recall definitions and facts which are important to us. This
material is widely known,  but quite rarely presented in a way which is most
convenient for our purposes. We concentrate the details that we use. For
similar introductions, we refer to \cite{alexaw}, \cite{beukers},
\cite{maintphd}.

\subsection{Algebraic curves}
\label{secalgcurves}

For general theory of algebraic curves we refer to \cite{fultonac} or to
\cite{shafarevich}. We assume algebraic curves to be reduced, irreducible,
smooth, complete (or projective), defined over $\CC$. In particular, the
projective line $\PP^1$ is $\CC\cup\{\infty\}$ set-theoretically. (All
curves in this paper are defined over $\QQ$, but we do not consider
airthmetic properties here.)

Let $C$ denote an algebraic curve. We denote the field of rational functions
on $C$ by $\CC(C)$. It can be generated by 2 functions, since $C$ is
birationally isomorphic to a (possibly singular) curve in $\PP^2$.
The function field $\CC(\PP^1)$ can be generated by 1 function; such a
generator is called a {\em rational parameter}.

If $P\in C$ and $f\in\CC(C)$, then $\mbox{ord}_P(f)$ denotes the {\em
valuation} of $f$ at $P$. If negative, this is the order of a pole of $f$ at
$P$; otherwise this is the vanishing order of $f$ at $P$. A {\em local
parameter} at $P$ is a function $t_P\in\CC(C)$ such that
$\mbox{ord}_P(t_P)=1$. For example, if
$C$ is the projective line, 
then $x-\alpha$ is a local parameter at the point $x=\alpha$, and $1/x$ is a
local parameter at $x=\infty$.

A {\em divisor} on $C$ is a finite formal sum $\sum_{P\in C} a_PP$, with
$a_P\in\ZZ$. The {\em degree} of such a divisor is the integer $\sum_{P\in
C} a_P$. The divisors form a commutative group under addition. For a
function $f\in\CC(C)$ we have its {\em principal divisor} $\sum_{P\in C}
\mbox{ord}_P(f)\,P$, which has degree $0$. Principal divisors form a
subgroup of degree 0 divisors. The quotient of these two groups is a {\em
Piccard group} of $C$; it is denoted by $\mbox{Pic}(C)$. For example,
$\mbox{Pic}(\PP^1)$ is the trivial group because all degree zero divisors on
$\PP^1$ are principal.

Explicit curves in this paper have either genus 0 (i.e., isomorphic to
$\PP^1$) or genus 1. Let $E$ denote a curve of genus 1. It can be
represented in a {\em Weierstrass form} $\xi^2=G_3(x)$, where $G_3(x)$ is a
cubic polynomial in $\CC[x]$. The point at infinity in this model by
$\mathcal{O}$. The Piccard group of $E$ is isomorphic (set-theoretically) to
$E$ itself. As usual, we identify a point $P\in E$ with the element of
$\mbox{Pic}(E)$ represented by the divisor $P-\mathcal{O}$. Then the
additive group law on $E$ can be given by the known chord-and-tangent
method. In particular, if three points of $E$ lie on one line of $\PP^2$,
they add up to the neutral element $\mathcal{O}$. The curve $E$ with this
group law is an {\em elliptic curve} $(E,\mathcal{O})$. Recall that a {\em
torsion point} on $E$ is a point of finite order. 
\begin{lemma} \label{princdiv}
Let $(E,\mathcal{O})$ denote an elliptic curve, 
and let $T=\sum_{P\in E} a_PP$ be a divisor on $E$. Then $T$ is a principal
divisor if and only if $\sum_{P\in E} a_PP=\mathcal{O}$ in the additive
group of $(E,\mathcal{O})$.
\end{lemma}
\proof Follows from the specified identification of $(E,\mathcal{O})$ with
$\mbox{Pic}(E)$. \qed

We also consider {\em radical functions} on $C$, that is, products of
$\bf{Q}$-powers of functions from $\CC(C)$. These are multi-valued
functions, but their branching points are poles or zeroes with finitely many
complex branches coming together. Valuations of those functions are well
defined at any point, and have values in $\QQ$ at the branching points.
Accordingly, we consider their principal divisors $\sum_{P\in C} a_PP$ with
coefficients $a_P\in\QQ$.
\begin{lemma} \label{principald}
Let $(E,\mathcal{O})$ denote an elliptic curve, 
and let $T=\sum_{P\in E} a_PP$ be a divisor with coefficients in
$\QQ$. 
Then $T$ is the principal divisor for a radical function if and only if the
following conditions hold:
\begin{enumerate}
\item The degree $\sum_{P\in E} a_P$ is zero. \item Let $\sum_{P\in E}
\widetilde{a}_PP$ be an integer multiple of $\,T$ such that all coefficients
$\widetilde{a}_P$ are integers. Then, in the additive group of $E$, the
point $\sum_{P\in C} \widetilde{a}_PP$ must be a torsion point.
\end{enumerate}
\end{lemma}
\proof Under these conditions, an integer factor $nT$ of $T$ would be a
principle divisor with integer coefficients. Then $T$ is a divisor of
$G^{1/n}$ for some $G\in\CC(E)$.

On the other hand, if $T$ is a divisor of a radical function $f$, then an
integer power $f^n$ is a rational function. The divisor $nT$ sums up to
$\mathcal{O}$ by Lemma \ref{princdiv}. Other integer factors of $T$ with
integral coefficients may sum up to a torsion point. \qed



\subsection{Finite coverings and pull-back transformations}
\label{apullbacks}

Consider a finite covering $\phi:C\to D$ from $C$ to other algebraic curve
$D$. It induces an algebraic field extension $\CC(C)\supset\CC(D)$. We
denote the degree of $\phi$ by $\deg\phi$. The genus $g(C)$ and $g(D)$ of
both curves and branching data are related by the Hurwitz formula:
\begin{equation}
2\,g(C)-2=\big( 2\,g(D)-2 \big)\,\deg\phi+\sum_{P\in C} \left( r_P-1
\right).\vspace{-6pt}
\end{equation}
Here $r_P$ is the branching order at $P$. It is equal to
$\mbox{ord}_P(t_{\phi(P)}\circ\phi)$.

Now we convene what we mean by a pull-back of hypergeometric equation
(\ref{hpgdee}) with respect to a finite covering. Let $C$ denote an
algebraic curve. Suppose that the function field $\CC(C)$ of $C$ is
generated by functions $x,\xi$. If $C$ is a rational curve, we may assume
that $\xi$ is not used and $x$ is a rational parameter of $C$.

Consider a finite covering $\phi:C\to\PP^1$. Let $z$ denote a rational
parameter for $\PP^1$. Then a pull-back of (\ref{hpgdee}) with respect to
$\phi$ is a differential equation is defined by transformation:
\begin{equation} \label{algtransf2}
z\longmapsto\phi(x,\xi), \qquad y(z)\longmapsto
Y(x,\xi)=\theta(x,\xi)\;y(\phi(x,\xi)).
\end{equation}
Here $\theta(x,\xi)$ is a radical function. Note that such a function has
the property that its logarithmic derivative $\theta'(x,\xi)/\theta(x,\xi)$
is in $\CC(C)$. We use the derivation on $\CC(C)$ that extends the usual
derivative on $\CC(x)$. If $C\cong\PP^1$, then transformation
(\ref{algtransf2}) is the following:
\begin{equation} \label{algtransf}
z\longmapsto\phi(x), \qquad y(z)\longmapsto Y(x)=\theta(x)\;y(\phi(x)).
\end{equation}

\subsection{Differential Galois theory}
\label{diffgalois}

A {\em differential field} $K$ is a field with a derivation, i.e., a map
$D:K\to K$ which satisfies $D(a+b)=D(a)+D(b)$ and the Leibnitz rule
$D(ab)=aD(b)+bD(a)$. One usually denotes $D(a)$ by
$a'$. 
An {\em extension} of the differential field $K$ is a differential field $L$
which contains $K$ and whose derivation extends the derivation of $K$. The
basic example of a differential field is the field $\CC(z)$ of rational
functions on $\PP^1$ with the usual derivation. Other example is the field
$\CC(C)$ of rational functions on an algebraic curve. To give a derivation
on $\CC(C)$ one may consider a finite covering $\phi:C\to\PP^1$ and the
corresponding unique extension of the usual derivation of $\CC(z)$.

Fix a differential field $K$ and consider a linear homogeneous differential
equation (\ref{genlhde}). Solutions in any extension of $K$ 
form a linear space over the {\em constant field} $\{a\in K\,|\,a'=0\}$. The
dimension of the solution space is at most $n$. A {\em Piccard-Vessiot
extension} $K_{PV}\supset K$ for (\ref{genlhde}) is, roughly speaking, a
minimal extension of differential fields, such that the solutions of
(\ref{genlhde}) in $K_{PV}$ form a linear space
of dimension $n$. 
The {\em differential Galois group} $G$ of (\ref{genlhde}) is the group of
autocoverings of $K$ that fix the elements of $K_{PV}$. The action of $G$ on
the $n$-dimensional space of solutions in $K_{PV}$ gives a faithful
$n$-dimensional representation of $G$. Therefore the differential Galois
group $G$ is usually considered as an algebraic subgroup of
$\mbox{GL}(n,\CC)$.

In \S \ref{darbouxcurves}, we utilize the {\em Riccati equation}
associated to (\ref{genlhde}). Solutions for the Riccati equation are
precisely the logarithmic derivatives $y'/y$ of solutions for
(\ref{genlhde}). Explicitly, the Riccati equation for (\ref{genlhde}) with
$n=2$ is $u'+u^2+a_1u+a_0=0$. Rational or algebraic solutions of the
Riccati equation are important in finding 
``closed form" solutions of the original equation (\ref{genlhde}), see
\cite{kovacic}. We refer to algebraic solutions of the Riccati equation in
our working definition of Darboux curves.

Suppose that hypergeometric equation (\ref{hpgdee}) has a finite monodromy $G$. 
Then the differential Galois group is isomorphic to $G$. (More
generally, the differential Galois group of a Fuchsian equation is
isomorphic to the Zariski closure of a representation of the monodromy
group.) The Piccard-Vessiot extension $K_{PV}\supset\CC(z)$ is a finite
Galois extension, the usual Galois group is isomorphic to $G$ as well. If
$y(z)\in K_{PV}$ is a solution of (\ref{hpgdee}), then $K_{PV}=\CC(z,y)$.

In most papers on differential Galois theory, 
second order differential equations are normalized to the form
$y(z)''=r(z)y(z)$, with $r(z)\in\CC(z)$. Hypergeometric equation
(\ref{hpgdee}) can be normalized by the transformation
\begin{equation} \label{normtr}
y(z)\longmapsto z^{(e_0-1)/2}\,(1-z)^{(e_1-1)/2}\,y(z).
\end{equation}
The normalized equation is:
\begin{equation} \label{normeq}
\frac{d^2y(z)}{dz^2}=\left( \frac{e_1^2-1}{4\,(z-1)^2}
+\frac{e_0^2-1}{4\,z^2}+\frac{1+e_\infty^2-e_0^2-e_1^2}{4\,z\,(z-1)} \right)
y(x).
\end{equation}
If the monodromy $G$ of a hypergeometric equation is isomorphic to
$A_4,S_4$ or $A_5$, then the differential Galois group of the normalized
equation (\ref{normeq}) is $G\times\{1,-1\}$. This does not change facts
that are important to us. Algebraic degree of Riccati solutions for
(\ref{hpgdee}) is the same as of Riccati solutions for (\ref{normeq}).


\subsection{Hypergeometric equations}
\label{hpgdeqs}

The singularities and local exponents of the hypergeometric equation  (\ref{hpgdee})
are conveniently revealed in Riemann's $P$-notation:
\begin{equation} \label{riemanp}
P\left\{ \begin{array}{ccc} 0 & 1 & \infty \\ 0 & 0 & a \\
1-c & c-a-b & b \end{array}\; z \; \right\}.
\end{equation}
As we see, the first row indicates the regular singular points, and the other
rows contain the local exponents and the variable $z$. 
The exponent differences obviously are $e_0=1-c$, $e_1=c-a-b$ and $e_\infty=a-b$.
These linear expressions in $a,b,c$ can be inverted, so that an ordered sequence
 of exponent differences determines the hypergeometric equation.
If we permute the exponent differences $e_0,e_1,e_\infty$ or
multiply some of them by $-1$, we get hypergeometric equations related by
well-known fractional-linear transformations \cite{specfaar}. In general,
there are 24 hypergeometric equations related in this way, and they share
the same (up to radical factors and fractional-linear change of the
independent variable) 24 hypergeometric Kummer's solutions.

A general basis of solutions for (\ref{hpgdee}) is
\begin{equation} \label{hsolsat0}
\hpg{2}{1}{a,\;b}{c}{\,z\,}, \qquad
z^{1-c}\;\hpg{2}{1}{a\!+\!1\!-\!c,\;b\!+\!1\!-\!c\,}{2-c}{\,z\,},
\end{equation}
where $\hpg{2}{1}{\!a,\,b}{c}{z}:=1+\frac{a\;b}{c\!\cdot 1!}\,z+
\frac{a(a+1)b(b+1)}{c(c+1)\cdot 2!}\,z^2 +\ldots$ is the Gauss
hypergeometric series. 


\subsection{Fuchsian equations}
\label{fuchsian}

All differential equations that we explicitly consider are {\em Fuchsian
equations}. 
These equations have only regular singular points. For equation
(\ref{genlhde}) this means the following: if $K=\CC(C)$ for an algebraic
curve $C$, then for any point $P\in C$ and for $i=1,\ldots,n$ we must have
$\mbox{ord}_P(a_i)\ge (n-i)\left(\mbox{ord}_P(t_P\!')-1\right)$, where $t_P$
is a local parameter at $P$. 
Local exponents at $P$ can be defined as follows: substitute $y=t_P^{\mu}$
into the Fuchsian equation and consider the terms to the power
$\mu+n(\left(\mbox{ord}_P(t_P\!')-1\right)$ of $t_P$ as an equation in
$\mu$; the roots of that equation are precisely the local exponents. The
local exponents at regular points are equal to $0,1,\ldots,n-1$.

In general, singularities and local exponents do not determine a Fuchsian
equation uniquely. Hence we cannot always use the $P$-notation for general
Fuchsian equations. However, in \S \ref{compevals} we write down arrays
of singularities and local exponents similar to (\ref{riemanp}).

\subsection{Contiguous relations of Gauss hypergeometric functions}
\label{contigrels}

Two Gauss hypergeometric functions 
are called {\em contiguous} (or {\em associated} in \cite{bateman}) if they
have the same argument $z$ and their parameters $a$, $b$ and $c$ differ
respectively by integers. As is known \cite[\S 2.5]{specfaar}, for any
three contiguous $\hpgo{2}{1}$ functions there is a {\em contiguous
relation}, which is a linear relation between the three functions where the
coefficients are rational functions in the parameters $a,b,c$ and the
argument $z$. A straightforward (though not efficient) method to compute a
contiguous expression for $\hpg{2}{1}{a+k,b+\ell}{c+m}{z}$ in terms of
$\hpg{2}{1}{a+1,b}{c}{z}$ and $\hpg{2}{1}{a,b}{c}{z}$ is the following. By
using the contiguous relations
\begin{eqnarray}
b\;\hpg{2}{1}{a,b+1}{c}{z}&=&
(b-a)\;\hpg{2}{1}{a,b}{c}{z}+a\;\hpg{2}{1}{a\!+1,b}{c}{z},\\
(c-1)\;\hpg{2}{1}{a,\,b}{c-1}{z}&=&(c-a-1)\;
\hpg{2}{1}{a,b}{c}{z}+a\;\hpg{2}{1}{a\!+1,b}{c}{z},
\end{eqnarray}
one eliminates the shifts in $b$ and $c$, and then by using the contiguous
relation
\begin{equation}
a(1-z)\,\hpg{2}{1}{\!a\!+1,b}{c}{z}=(2a-c-az+bz)\;\hpg{2}{1}{\!a,b}{c}{z}
+(c-a)\;\hpg{2}{1}{\!a\!-1,b}{c}{z}
\end{equation}
one gets an expression with two contiguous terms. Effective computation of
contiguous relations is considered in \cite{contigrels}. They can be
computed in $O(\log(\max(k,l,m))$ steps, but complexity of expressions in
each such step grows exponentially, and the output is $O(\max(k,l,m))$.

One can rewrite contiguity conditions in terms of local exponent differences
at $z=0$, $1$, $\infty$ for the hypergeometric equation, since the
parameters $a$, $b$, $c$ determine the exponent differences and vice
versa (if the sign of exponent differences is taken into account). The
main hypergeometric solutions (\ref{hsolsat0}) of two hypergeometric
equations (\ref{hpgdee}) are contiguous if for each $X\in\{0,1,\infty\}$ the
difference of signed exponent differences at $X$ of the two equations
is an integer, and the sum of the three integer differences is even. Two
hypergeometric equations have solutions contiguous to each other (or
equivalently, they have the same Schwarz type) if one can choose a
permutation of exponent differences and their sign in such a way that
the just described situation occurs.

For example, the parameters $e_0,e_1,e_\infty$ of hypergeometric equations
of the Schwarz type $(1/3,1/3,2/3)$ can be characterized as follows: they
are rational numbers, their denominators are equal to 3, and the sum of
their numerators is even.

\subsection{Fiber products of curves}
\label{fiberprod}


Let $C_1$ and $C_2$ denote two curves over $\CC$. Let $\phi_1:C_1\to\PP^1$
and $\phi_2:C_2\to\PP^1$ be two finite coverings of degree $m$ and $n$
respectively.
The fiber product of $\phi_1:C_1\to\PP^1$ and $\phi_2:C_2\to\PP^1$ is a
curve $B$ with two coverings $\psi_1:B\to C_1$ and $\psi_2:B\to C_2$ such
that $\phi_1\circ\psi_1=\phi_2\circ\psi_2$, and for any other curve
$\widetilde{B}$ with coverings $\widetilde{\psi}_1:\widetilde{B}\to C_1$ and
$\widetilde{\psi}_2:\widetilde{B}\to C_2$ satisfying
$\phi_1\circ\widetilde{\psi}_1=\phi_2\circ\widetilde{\psi}_2$ there is a
unique covering $\xi:\widetilde{B}\to B$ such that
$\widetilde{\psi}_1=\psi_1\circ\xi$ and $\widetilde{\psi}_2=\psi_2\circ\xi$.
Then the following diagram commutes:
\[
\begin{picture}(270,100)
\put(145,-1){$\PP^1$} \put(115,40){\vector(1,-1){30}}
\put(185,40){\vector(-1,-1){30}} \put(106,45){$C_1$} \put(187,45){$C_2$}
\put(145,86){\vector(-1,-1){30}} \put(157,86){\vector(1,-1){30}}
\put(147,90){$B$} \put(71,90){$\widetilde{B}$}
\put(75,86){\vector(1,-1){30}} \put(81,87){\vector(3,-1){102}}
\put(82,93){\vector(1,0){61}} \put(119,21){$\phi_1$}\put(174,21){$\phi_2$}
\put(111,95){$\zeta$} \put(80,62){$\widetilde{\psi}_1$}
\put(140,75){$\psi_1$} \put(152,50){$\widetilde{\psi}_2$}
\put(172,75){$\psi_2$}
\end{picture}
\]
We have $\deg\psi_1=\deg\phi_2$ and $\deg\psi_2=\deg\phi_1$. 
On the level for function fields, we have
$\CC(B)=\CC(C_1)\otimes_{\scriptsize{\bf C}(\PP^1)}\CC(C_2)$.

A birational model for $B$ is the curve on $C_1\times C_2$ of those points
$(X_1,X_2)\in C_1\times C_2$ which satisfy $\varphi_1(X_1)=\varphi_2(X_2)$.
This is a singular model in general. A singular point corresponds to a pair
$(X_1,X_2)\in B$ such that $X_1$ and $X_2$ have branching indices $r_1>1$,
$r_2>1$ respectively (with respect to $\phi_1$ and $\phi_2$). Such a
singularity is of type $x_1^{r_1}-x_2^{r_2}$; by resolving it we get
$\gcd(r_1,r_2)$ points that correspond to $(X_1,X_2)$ on a non-singular
model for $B$. 
This information allows us to compute the branching data for the projections
$\psi:B\to C_1$ and $\psi:B\to C_2$ and the genus
of $B$. 
For example, if $X_1\in C_1$ has branching index $r_1$ (with respect to
$\phi_1$), and the branching data of $\phi_2$ above $\phi_1(X_1)$ is
$a_1+\ldots+a_k$, then the branching data for $\psi_1$ above $X_1$ is the
following: $\gcd(a_1,r_1)*\frac{\mbox{lcm}(a_1,r_1)}{r_1}+
\ldots+\gcd(a_k,r_1)*\frac{\mbox{lcm}(a_k,r_1)}{r_1}$.

\small

\bibliographystyle{alpha}
\bibliography{../../hypergeometric}

\end{document}